\documentclass[10pt]{article}

\usepackage{dirtytalk}

\usepackage{amsfonts}
\usepackage{amssymb}
\usepackage{amsmath}
\usepackage{enumerate}
\usepackage{amsthm}

\newcommand{\cf}{\mathrm{cf}}

\newcommand{\cov}{\mathrm{cov}}

\newcommand{\mad}{\mathrm{MAD}}
\newcommand{\pp}{\mathrm{pp}}

\newtheorem{Th}{\bf THEOREM}[section]

\newtheorem{Lem}[Th]{\bf LEMMA}

\newtheorem{Pro}[Th]{\bf PROPOSITION}

\newtheorem{fact}[Th]{\bf FACT}

\newtheorem{Cor}[Th]{\bf COROLLARY}

\newtheorem{Obs}[Th]{\bf OBSERVATION}

\theoremstyle{definition} \newtheorem{Def}[Th]{\bf DEFINITION}
 
\theoremstyle{remark}\newtheorem{Rmk}[Th]{\bf REMARK}

\theoremstyle{question}\newtheorem{Ques}[Th]{\bf QUESTION}

\newtheorem{Conj}[Th]{\bf CONJECTURE}

\usepackage[english]{babel}
\makeindex

\title{TOWERS AND CLUBS}
\author{Pierre MATET}

\date{}

\begin{document}

\maketitle

\renewcommand{\thefootnote}{\arabic{footnote}} 	

\renewcommand{\thefootnote}{}                                



\vskip 0,7cm

\begin{abstract}  We revisit several results concerning club principles and nonsaturation of the nonstationary ideal, attempting to improve them in various ways. So we typically deal with a (non necessarily normal) ideal $J$ extending the nonstationary ideal on a regular uncountable (non necessarily successor) cardinal $\kappa$, our goal being to witness the nonsaturation of $J$ by the existence of towers (of length possibly greater than $\kappa^+$).     
 \end{abstract}
\bigskip

\section{INTRODUCTION}

\medskip

We will show that by modifying the proofs of some well-known results concerning non-saturation of the nonstationary ideal $NS_\kappa$, one may obtain towers in $(P (\kappa) / NS_\kappa, \subseteq)$ or/and $(P (\kappa) / NS_\kappa, \supseteq)$. Since these proofs usually involve one form or another of Club, we are led to revisit a number of results concerning this principle and its (many) variants.  
\bigskip

\section{DIAMOND LITE}
\bigskip

\subsection{Ideals and density}

\medskip 

Let us first recall some definitions and facts that will be needed later. We start with ideals.

\medskip

\begin{Def} For a set $A$ and a cardinal $\rho$, we set $P_\rho (A) = \{ a\subseteq A :  \vert a \vert < \rho\}$ and $ [A]^\rho = \{ a\subseteq A :  \vert a \vert  = \rho\}$.
\end{Def}

{\bf Throughout the paper $\kappa$ will denote a regular uncountable cardinal.}

\begin{Def} By an {\it ideal  on} $\kappa$ we mean a nonempty collection $J$ of subsets of $\kappa$ such that 
 \begin{itemize}
\item $\kappa \notin J$.
\item $\kappa \subseteq J$ ;
\item $P(A)\subseteq J$ for all $A\in J$.  
\item $A\cup B\in J$ whenever $A, B\in J$. 
\end{itemize}

Given an ideal $J$ on $\kappa$, we let $J^+ = P(\kappa) \setminus J$, $J^\ast = \{ A\subseteq \kappa : \kappa \setminus A\in J\}$, and $J \vert A = \{ B\subseteq \kappa : B \cap A\in J\}$  for each $A\in J^+$.  $J$ is {\it prime} if $J^+ = J^\ast$, and {\it nowhere prime} if $J \vert A$ is prime for no $A \in J^+$. 
$J$ is {\it $\kappa$-complete} if  $\bigcup Q \in J$ for every $Q \in P_\kappa (J)$. For a cardinal $\rho$ and $Y\subseteq P(\kappa)$, $J$ is {\it $Y$-$\rho$-saturated} if there is no  $Q \subseteq J^+$ with $\vert Q \vert = \rho$ such that $A \cap B \in Y$ for any two distinct members $A, B$ of $Q$.
$J$ is {\it $\rho$-saturated} if it is $J$-$\rho$-saturated. 

$J$ is {\it subnormal} if $J \subseteq K$ for some normal ideal $K$ on $\kappa$.
 \end{Def}

\begin{Def}  We let $I_\kappa$ and $NS_\kappa$ denote, respectively, the noncofinal ideal on $\kappa$ and  the nonstationary ideal on $\kappa$.

We let ${\cal C}_\kappa$ denote the collection of all closed unbounded subsets of $\kappa$.

For $A \subseteq \kappa$, we let $acc (A) = \{\alpha \in \kappa \setminus \{0\} : \sup (A \cap \alpha) = \alpha\}$.
\end{Def}

\begin{Def} Given a cardinal $\theta$, $E_\theta^\kappa$ (respectively, $E_{< \theta}^\kappa$, $E_{\geq \theta}^\kappa$) denotes the set of all $\alpha \in acc (\kappa)$ with $\cf (\alpha) = \theta$ (respectively, $\cf (\alpha) < \theta$, $\cf (\alpha) \geq \theta$).
\end{Def}

\begin{Def}  Let $S$ be a stationary subset of $\kappa$. For $\gamma \in E^\kappa_{\geq \omega_1}$, $S$ {\it reflects at} $\gamma$ if $S \cap \gamma$ is stationary in $\gamma$.
 \end{Def}
 
\begin{Rmk} If $S \subseteq E^\kappa_\theta$ reflects at $\gamma$, then $\cf (\gamma) > \theta$.
\end{Rmk}

\medskip

We next turn to density numbers and meeting numbers.

\medskip

\begin{Def}  Given two cardinals  $\tau\leq\sigma$ with $1 \leq \tau$ and $\omega \leq \sigma$, $d(\tau,\sigma)$ (respectively, $m(\tau, \sigma)$) denotes the least cardinality of any $X\subseteq [\sigma]^\tau$  with the property that  for any $e\in [\sigma]^\tau$,  there is $x\in X$  with $x\subseteq e$ (respectively, $\vert x \cap e \vert = \tau$).
\end{Def}

\begin{Rmk} Thus $d(\tau,\sigma) =$ the cofinality of the poset $([\sigma]^\tau, \supseteq)$.
\end{Rmk}

\begin{Def} 
 Given two infinite cardinals $\tau\leq\sigma$, $\mad_{\tau,\sigma}$ denotes the collection of all $Q\subseteq [\sigma]^\tau$  such that  
 \begin{itemize}
 \item $\vert a\cap b \vert < \tau$ for any two distinct members $a,b$ of $Q$ ;
 \item  for any $c\in [\sigma]^\tau$,  there is $a\in Q$  with $\vert a\cap c \vert = \tau$.
\end{itemize} 
\end{Def}
 
\begin{fact}
{\rm (\cite{Koj}, \cite{Diamondstar}, \cite{GMS})}   Let $\tau \leq \sigma$  be two infinite cardinals. Then the following hold~: 
\begin{enumerate} [\rm (i)]
\item $\sigma \leq d(\tau, \sigma) \leq \sigma^\tau$. 
\item $d(\tau,\sigma) = \sigma$  if and only if $\cf (\sigma) \not= \cf (\tau)$ and $d(\tau, \chi)\leq \sigma$  for any cardinal $\chi$  with $\tau\leq\chi < \sigma$.
\item Suppose that $\chi^{<\tau}\leq\sigma < \chi^\tau$  for some cardinal $\chi$. Then $d(\tau,\sigma) = \sigma^\tau$.
\item $d(\tau, \chi) \leq d(\tau, \sigma)$ for any cardinal $\chi$ with $\tau \leq \chi \leq \sigma$.
\item $d(\tau,\sigma^+) = \max \{d(\tau,\sigma), \sigma^+\}$.
\item Suppose that $\sigma$ is a limit cardinal with $\cf (\sigma) \not= \cf (\tau)$. Then $d(\tau,\sigma) = \sup \{ d(\tau, \chi) : \tau \leq \chi < \sigma\}$.
 \item If $\sigma < \tau^{+\cf (\tau)}$, then $d(\tau,\sigma) = \max \{d(\tau,\tau),\sigma\}$.
 \item $d(\tau,\sigma) \geq \vert Q \vert$  for all $Q\in \mad_{\tau,\sigma}$.
 \item If $d(\tau,\tau) < d(\tau,\sigma)$,  then $\vert Q \vert = d(\tau,\sigma)$  for all $Q\in \mad_{\tau,\sigma}$.
 \item $d(\tau, \sigma) = \max \{ d(\tau, \tau), m(\tau, \sigma)\}$.
\end{enumerate}
\end{fact} 

\begin{Def}  {\it Shelah's Strong Hypothesis} {\rm (SSH)} asserts that $\pp(\chi) = \chi^+$ for every singular cardinal $\chi$.
\end{Def}

\begin{fact} {\rm (\cite{SSH})}  The following are equivalent :
 \begin{enumerate}[\rm (i)]
 \item Shelah's Strong Hypothesis.
 \item Given two infinite cardinals $\tau < \sigma$, $m(\tau,\sigma)$ equals $\sigma$ if $\cf (\sigma) \not= \cf (\tau)$, and $\sigma^+$ otherwise.
 \end{enumerate}
 \end{fact}

\begin{fact} {\rm (\cite{SheRGCH}, \cite{Diamondstar})} \begin{enumerate}[(i)] 
\item Let $\rho$  be an uncountable strong limit cardinal, and $\sigma\geq\rho$  be a cardinal. Then there is $\alpha < \rho$  such that for any infinite cardinal $\tau$  with $\alpha\leq\tau\leq\rho$, $d(\tau,\sigma)$ equals $\sigma$ if $\alpha\leq\cf (\tau)$, and $d(\cf (\tau),\sigma)$  otherwise.
\item  Suppose that $\rho < \kappa$  is an uncountable strong limit cardinal, and $\kappa$ is a limit cardinal. Then we may find $\chi < \rho$  with the following property : If $\tau$  is a regular cardinal with $\chi\leq\tau < \rho$,  then $d(\tau,\sigma) < \kappa$  for every cardinal $\sigma$ with $\tau\leq \sigma < \kappa$.
\end{enumerate}
\end{fact}

\subsection{J'enl\`eve le haut}

\medskip

Our starting point is a result of Gregory on diamond star. The following guessing principles were introduced by Jensen \cite{Jen}.

\medskip

\begin{Def}  Given a $\kappa$-complete ideal $J$ on $\kappa$, $\diamondsuit_\kappa^\ast [J]$ (respectively, $\diamondsuit_\kappa^- [J]$) asserts the existence of $t^{i}_\alpha \subseteq \alpha$ for $i < \alpha < \kappa$ such that $\{\alpha < \kappa : \exists i < \alpha (t^{i}_\alpha = A \cap \alpha)\}$ lies in $J^\ast$ (respectively, $J^+$) for every $A \subseteq \kappa$.

$\diamondsuit_\kappa [J]$ asserts the existence of $s_\alpha \subseteq \alpha$ for $\alpha < \kappa$ such that $\{\alpha < \kappa : s_\alpha = A \cap \alpha\}$ lies in $J^+$ for every $A \subseteq \kappa$.
\end{Def}

\begin{Rmk} \begin{enumerate}[(i)] 
\item If $\diamondsuit_\kappa^\ast [J]$ holds, then so does $\diamondsuit_\kappa^\ast [K]$ for any $\kappa$-complete ideal $K$ on $\kappa$ extending $J$.
\item $\diamondsuit_\kappa^\ast [J] \Rightarrow \diamondsuit_\kappa^- [J]$.
\item $\diamondsuit_\kappa [J] \Rightarrow \diamondsuit_\kappa^- [J]$.
\end{enumerate} 
\end{Rmk}

\begin{fact} \begin{enumerate}[(i)] 
\item {\rm (\cite{Kun}, \cite{Shelah1})} Suppose that either $\kappa$ is a successor cardinal, or $J$ is normal. Then $\diamondsuit_\kappa^- [J]  \Rightarrow \diamondsuit_\kappa [J]$.
\item  {\rm (\cite{Jen}, \cite{Shelah1})} $\diamondsuit_\kappa^- [J] \Rightarrow 2^{< \kappa} = \kappa$.
\item  {\rm (Folklore)} If $\diamondsuit_\kappa [J]$ holds, then $J$ is not $I_\kappa$-$2^\kappa$-saturated.
\end{enumerate} 
\end{fact}

\medskip

Gregory's result \cite{Gre} asserted that if $\kappa = \nu^+ = 2^\nu$, then $\diamondsuit^\ast_\kappa [NS_\kappa \vert E^\kappa_\theta]$ holds for any regular infinite cardinal $\theta$ with $\nu^\theta = \nu$. It was later strengthened by Shelah (\cite{She78}, \cite{She81}) and others (\cite{Rin11}, \cite{Diamondstar}). Its present form (not necessarily the final one) reads as follows.

\medskip

\begin{fact} Suppose that $\kappa = \nu^+ = 2^\nu$, and let $\theta$ be a regular infinite cardinal less than $\nu$ such that $d (\theta, \nu) = \nu$. Then $\diamondsuit^\ast_\kappa [NS_\kappa \vert E^\kappa_\theta]$ holds.  
\end{fact}

\medskip

Let us discuss the requirement that $d (\theta, \nu) = \nu$. By Facts 2.10 (x) and  2.12, under SSH, it reduces to the condition that $d (\theta, \theta) \leq \nu$ (which will be satisfied if $\kappa$ is large enough) and $\theta \not= \cf (\nu)$. On the other hand, if there is a strong limit cardinal $\tau$ with $\theta < \tau < \kappa$, and $\theta$ is large enough, then by Fact 2.13, $d (\theta, \nu) = \nu$ will hold. So there are many cases when the condition $d (\theta, \nu) = \nu$ is satisfied. But what can be said when it is not ? Shelah has the following answer.

\begin{fact}  {\rm(\cite{She10})} Suppose that $\kappa = \nu^+ = 2^\nu$, and let $\theta$ be a regular infinite cardinal less than $\kappa$ with $\theta \not= \cf (\nu)$. Then $\diamondsuit_\kappa [J]$ holds for any $\kappa$-complete ideal on $\kappa$ extending $NS_\kappa \vert E^\kappa_\theta$.  
\end{fact}

\medskip

Notice that the result also applies to ideals that are not normal. Of course if $\diamondsuit_\kappa [K]$ holds for some normal ideal on $\kappa$, then so does $\diamondsuit_\kappa [J]$ for any $\kappa$-complete ideal $J$ on $\kappa$ included in $K$. So the ideals that would not be covered if the result were only stated for normal ideals are those that are not subnormal. The following result provides a description of these ideals.

\medskip

\begin{fact}  {\rm(\cite{BTW})} Given a $\kappa$-complete ideal $J$ on $\kappa$, the following are equivalent :
\begin{enumerate}[\rm (i)]
\item $J$ is not subnormal.  
\item There is a partition $\langle S_\alpha : \alpha < \kappa \rangle$ of $\kappa \setminus \{0\}$ into stationary sets $S_\alpha$ with $S_\alpha \cap (\alpha + 1) = \emptyset$ such that $J$ extends the $\kappa$-complete ideal generated by $NS_\kappa \cup \{S_\alpha : \alpha < \kappa\}$.
\end{enumerate}
\end{fact}

\medskip

Let us return to Shelah's result. How does it look like if we go further and remove the remaining cardinal arithmetic hypothesis ? This paper originated in our desire to prove the following.

\medskip

\begin{Conj}  Suppose that $\kappa = \nu^+$, and let $\theta$ be a regular infinite cardinal less than $\kappa$ with $\theta \not= \cf (\nu)$. Then no $\kappa$-complete ideal on $\kappa$ extending $NS_\kappa \vert E^\kappa_\theta$ is $I_\kappa$-$\kappa^+$-saturated.
 \end{Conj}

\medskip

Why only $\kappa^+$ ? Just to play it safe, since $2^\kappa$ would not be suitable (Foreman and Magidor \cite{FM} showed that if $V = L$ and $\sigma$ Cohen subsets of $\omega_1$ are added, where $\sigma \geq \kappa^{++}$, then in the extension, $NS_\kappa$ is $\kappa^{++}$-saturated).
\bigskip

\subsection{J'enl\`eve le bas}


\medskip

\begin{Def}  Given a $\kappa$-complete ideal $J$ on $\kappa$, $\clubsuit_\kappa [J]$) asserts the existence of  $s_\alpha \subseteq \alpha$ with $\sup s_\alpha = \alpha$ for $\alpha \in acc (\kappa)$ such that $\{ \alpha \in acc (\kappa) : s_\alpha \subseteq A\} \in J^+$ for all $A \in [\kappa]^\kappa$.

The principle $\clubsuit_\kappa^\ast [J]$ asserts the existence of $s^{i}_\delta \subseteq \alpha$ with $\sup B^{i}_\alpha = \alpha$  for $i \in \alpha \in acc ( \kappa)$ such that $\{ \alpha < \kappa : \exists i < \alpha (s^{i}_\alpha \subseteq A)\} \in J^\ast$ for all $A \in [\kappa]^\kappa$.
\end{Def}

\medskip

$\clubsuit_{\omega_1} [NS_{\omega_1}]$ is usually denoted by $\clubsuit$ and known as Ostaszewski's guessing principle. 

It is easy to see that if $J$ extends $NS_\kappa$, then $\diamondsuit_\kappa [J]$ (respectively, $\diamondsuit^\ast_\kappa [J]$) implies $\clubsuit_\kappa [J]$ (respectively, $\clubsuit_\kappa^\ast [J])$. By a result of Devlin (see \cite{Ost}), $\diamondsuit_{\omega_1} [NS_{\omega_1}]$ follows from CH + $\clubsuit$. This easily generalizes.

\medskip

\begin{Obs} Given a $\kappa$-complete ideal $J$ on $\kappa$ extending $NS_\kappa$, the following are equivalent :
 \begin{enumerate}[\rm (i)]
 \item $\diamondsuit_\kappa [J]$ holds.
\item $\clubsuit_\kappa [J]$ holds and $2^{< \kappa} = \kappa$.
\end{enumerate}   
 \end{Obs}

{\bf Proof.}  By the proof of Observation 3.4 below.
 \hfill$\square$

\medskip

The starred version is established by a similar argument.

\medskip

 \begin{Obs} Given a $\kappa$-complete ideal $J$ on $\kappa$ extending $NS_\kappa$, the following are equivalent :
 \begin{enumerate}[\rm (i)]
 \item $\diamondsuit^\ast_\kappa [J]$ holds.
\item $\clubsuit^\ast_\kappa [J]$ holds and $2^{< \kappa} = \kappa$.
\end{enumerate}   
 \end{Obs} 

\medskip

Observation 2.23 improves a result of \cite{Fuchs} that asserts that if $\kappa = \nu^+$ and $\tau$ is a regular infinite cardinal less than $\nu$ such that $d(\tau, \sigma) \leq \nu$ for every cardinal $\sigma$ with $\tau \leq \sigma < \nu$, then for any $S \in NS^+_\kappa \cap P(E^\kappa_\tau)$, $\diamondsuit^\ast_\kappa [NS_\kappa \vert S]$ holds just in case $\clubsuit^\ast_\kappa [NS_\kappa \vert S]$ holds and $2^{< \kappa} = \kappa$.

\medskip

The consistency of $\clubsuit$ with the negation of the Continuum Hypothesis (and therefore with the negation of $\diamondsuit_{\omega_1} [NS_{\omega_1}]$) has been established by Shelah \cite{She80}.

\begin{Def}  Given a $\kappa$-complete ideal $J$ on $\kappa$, $\clubsuit_\kappa^{\rm ev} [J]$ asserts the existence of  $s_\alpha \subseteq \alpha$ with $\sup s_\alpha = \alpha$ for $\alpha \in acc (\kappa)$ such that $\{ \alpha \in acc (\kappa) : \exists \beta < \alpha (s_\alpha \setminus \beta \subseteq A)\} \in J^+$ for all $A \in [\kappa]^\kappa$.
\end{Def}

Obviously, $\clubsuit_\kappa [J] \Rightarrow \clubsuit_\kappa^{\rm ev} [J]$. The principle $\clubsuit_{\omega_1}^{\rm ev} [NS_{\omega_1}]$ is denoted by $\clubsuit_w$ in \cite{FSS}, and by $\clubsuit^1$ in \cite{DZ2} where its consistency with the negation of $\clubsuit$ is established.

It is known (see \cite{DZ1}, \cite{FSS}) that for any $S \in NS^+_\kappa$, $\clubsuit_\kappa [NS^+_\kappa \vert S]$ holds if and only if there is $s_\alpha \subseteq \alpha$ with $\sup s_\alpha = \alpha$ for $\alpha \in acc (\kappa)$ such that $\{ \alpha \in S \cap acc (\kappa) : s_\alpha \subseteq A)\} \not= \emptyset$ for all $A \in [\kappa]^\kappa$. This works also for the eventual-guessing variant.
\medskip

\begin{Obs} Given $S \in NS^+_\kappa$, the following are equivalent :
 \begin{enumerate}[\rm (i)]
 \item $\clubsuit_\kappa^{\rm ev} [NS^+_\kappa \vert S]$.
\item There is $s_\alpha \subseteq \alpha$ with $\sup s_\alpha = \alpha$ for $\alpha \in acc (\kappa)$ such that for any $A \in [\kappa]^\kappa$, $\{ \alpha \in S \cap acc (\kappa) : \exists \beta < \alpha (s_\alpha \setminus \beta \subseteq A)\} \not= \emptyset$.
\end{enumerate}   
 \end{Obs}

{\bf Proof.}  (i) $\rightarrow$ (ii) : Trivial.

(ii)$ \rightarrow$ (i) : Let $\langle s_\alpha : \alpha \in acc (\kappa) \rangle$ be as in (ii), and suppose toward a contradiction that (i) fails. Then there must be $A \in [\kappa]^\kappa$ and $C \in {\cal C}_\kappa \cap P (acc (\kappa))$ such that $(s_\alpha \setminus \beta) \setminus A \not= \emptyset$ whenever $\alpha \in C \cap S$ and $\beta < \alpha$. Set $B = \{ \min (A \setminus \alpha) : \alpha \in C\}$. We may find $\alpha \in S \cap acc (\kappa)$ and $\beta < \alpha$ such that $s_\alpha \setminus \beta \subseteq B \subseteq A$. But then $C$ is cofinal in $\alpha$, and consequently $\alpha \in C$. Contradiction.
 \hfill$\square$ 

\medskip

\begin{Rmk} Assuming GCH in $V$, Baumgartner \cite{Baum} has constructed a cofinality-preserving generic extension in which there is $S_i \in [\kappa]^\kappa$ for $i < \kappa^+$ such that $\vert S_i \cap S_j \vert < \aleph_0$ whenever $i < j < \kappa^+$. In this extension, $\clubsuit_\kappa^{\rm ev} [NS_\kappa \vert E^\kappa_\theta]$ must fail for every regular infinite cardinal $\theta < \kappa$, since it is well-known (see e.g.\cite{FSS}) that $\clubsuit_\kappa^{\rm ev} [NS_\kappa \vert E^\kappa_\theta]$ implies the existence of $X \subseteq [\kappa]^\theta$ with $\vert X \vert = \kappa$ such that for any $A \in [\kappa]^\kappa$, there is $x \in X$ with $x \subseteq A$. \end{Rmk}

\medskip

Garti \cite{Gar} observed that it follows from $\clubsuit$ that $NS_{\omega_1}$ is not $I_{\omega_1}$-$\omega_2$-saturated. An easy modification of his proof yields the following.\medskip

\begin{Obs} Let $J$ be a $\kappa$-complete ideal on $\kappa$ such that $\clubsuit_\kappa^{\rm ev} [J]$ holds, and $\rho$ be an infinite cardinal such that $I_\kappa$ is not $\rho$-saturated. Then $J$ is not $I_\kappa$-$\rho$-saturated.
 \end{Obs}

{\bf Proof.} Let $s_\alpha \subseteq \alpha$ with $\sup s_\alpha = \alpha$ for $\alpha \in acc (\kappa)$ be such that for any $A \in [\kappa]^\kappa$, $S_A \in J^+$, where 

\centerline{$S_A = \{ \alpha \in acc (\kappa) : \exists \beta < \alpha (s_\alpha \setminus \beta \subseteq A)\}$.}

 Pick $A_i \in [\kappa]^\kappa$ for $i < \rho$ so that $\vert A_i \cap A_j \vert < \kappa$ whenever $i < j < \rho$. We claim that $\vert S_{A_i} \cap S_{A_j} \vert < \kappa$ whenever $i < j < \rho$. Suppose otherwise, and fix $i < j < \kappa$ such that $\vert S_{A_i} \cap S_{A_j} \vert = \kappa$. Inductively define $\alpha_\xi \in S_{A_i} \cap S_{A_j}$ and $\gamma_\xi \in s_{\alpha_\xi} \cap A_i \cap A_j$ for $\xi < \kappa$ so that $\gamma_\xi > \sup \{ \alpha_\eta : \eta < \xi \}$. Then $\{ \gamma_\xi : \xi < \kappa \}$ is a size $\kappa$ subset of $A_i \cap A_j$. Contradiction.
 \hfill$\square$ 

\medskip

Let us now introduce the kind of towers we will be working with.

\medskip

\begin{Def} Given an ideal $J$ on $\kappa$,  $Y\subseteq P(\kappa)$ and an ordinal $\delta$, a {\it descending} (respectively, {\it ascending}) $(J, Y)${\it -tower of length} $\delta$ is a sequence $\langle A_\alpha : \alpha < \delta \rangle$ such that
\begin{itemize}
\item $A_\alpha \in J^+$ for all $\alpha < \delta$ ;
\item $A_\beta \setminus A_\alpha \in Y$ (respectively, $A_\alpha \setminus A_\beta \in Y$) and $A_\alpha \setminus A_\beta \in J^+$ (respectively,  $A_\beta \setminus A_\alpha \in J^+$) whenever $\alpha < \beta < \delta$.
\end{itemize}

A descending (respectively, ascending) $(J, Y)$-tower is {\it maximal} if there is no  descending (respectively, ascending) $(J, Y)$-tower properly extending it.
\end{Def}

\begin{Obs}  Let $\tau \geq 2$ be a cardinal such that there exists a descending (respectively ascending) $(J, Y)$-tower of length $\tau$. Then $J$ is not  $Y$-$\tau$-saturated.
\end{Obs}

{\bf Proof.} Let $\langle A_\alpha : \alpha < \tau \rangle$ be an ascending $(J, Y)$-tower. For $\alpha < \tau$, set $S_\alpha = A_{\alpha + 1} \setminus A_\alpha$. Then clearly $\{S_\alpha : \alpha < \tau\} \subseteq J^+$. Furthermore $S_\gamma \cap S_\alpha \subseteq A_{\gamma + 1} \setminus A_\alpha$ whenever $\gamma < \alpha < \tau$. Descending towers are handled in a similar way.
\hfill$\square$  

\begin{Def}  We let $Depth ([\kappa]^\kappa, \nearrow)$ (respectively $Depth ([\kappa]^\kappa, \searrow)$) denote the least ordinal $\eta$ such that there is no ascending (respectively, descending) $(I_{\kappa}, I_{\kappa})$-tower of length $\eta$.
\end{Def}

\begin{Th}  Let $J$ be a $\kappa$-complete ideal on $\kappa$ extending $NS_\kappa$, and $\eta$ be an infinite ordinal less than $Depth ([\kappa]^\kappa, \searrow)$ (respectively $Depth ([\kappa]^\kappa, \nearrow)$). Suppose that $\clubsuit^{\rm ev}_\kappa [J]$ holds. Then there exists a descending (respectively, ascending) $(J, I_ \kappa)$-tower of length $\eta$.
\end{Th}

{\bf Proof.} Select $s_\alpha \subseteq \alpha$ with $\sup s_\alpha = \alpha$ for $\alpha \in acc (\kappa)$ such that 

\centerline{$\{\alpha : \exists \beta < \alpha (s_\alpha \setminus \beta \subseteq A)\} \in J^+$}

for all $A \in [\kappa]^\kappa$. For $A \in [\kappa]^\kappa$, let $S_A$ denote the set of all $\alpha \in acc (\kappa)$ such that $\sup (A \cap \alpha) = \alpha > \sup (s_\alpha \setminus A)$. Notice that $S_A \in J^+$. Now let $\langle A_i : i < \eta \rangle$ be an ascending $(I_\kappa, I_\kappa)$-tower. Fix $i < j < \eta$, and pick $\delta < \kappa$ so that $A_i \setminus \delta \subseteq A_j$. Then $S_{A_i} \setminus (\delta + 1) \subseteq S_{A_j}$, and consequently $\vert S_{A_i} \setminus S_{A_j} \vert < \kappa$. Moreover, $S_{A_j \setminus A_i} \subseteq S_{A_j} \setminus S_{A_i}$. Thus $\langle S_{A_i} : i < \delta \rangle$ is an ascending $(J, I_\kappa)$-tower. The descending case is left to the reader. 
\hfill$\square$

\medskip

Thus $\clubsuit^{\rm ev}_\kappa [J]$ transmutes almost disjoint families of subsets of $\kappa$ into almost disjoint families of sets in $J^+$ of the same power, and descending (respectively, ascending) $(I_{\kappa}, I_{\kappa})$-towers into descending (respectively, ascending) $(J, I_ \kappa)$-towers of the same length. 

\bigskip

\subsection{Depths}

\medskip

To give some background to Theorem 2.31, in this subsection we discuss the existence of ascending (respectively, descending) towers.

\medskip

\begin{Def}  Given $f, g \in {}^\kappa \kappa$ , $f <^\ast g$ means that 

\centerline{$\vert \{\alpha < \kappa : f (\alpha) \geq g (\alpha\} \vert < \kappa$.} 
\end{Def}

\begin{Def}  We let $\frak{b}_\kappa$ denote the least cardinality of any $F \subseteq {}^\kappa \kappa$ with the property that there is no $g \in {}^\kappa \kappa$ such that $f <^\ast g$ for all $f \in F$. 
\end{Def}

\begin{Rmk} By an argument that goes back to Rothberger, there exist $f_\alpha \in  {}^\kappa \kappa$ for $\alpha < \frak{b}_\kappa$ such that 
\begin{itemize}
\item $f_\alpha <^\ast f_\beta$ for $\alpha < \beta < \frak{b}_\kappa$ ; 
\item there is no $g \in {}^\kappa \kappa$ such that $f_\alpha <^\ast g$ for all $\alpha \in \frak{b}_\kappa$. 
\end{itemize}
Notice that it follows that $\frak{b}_\kappa$ is regular.
\end{Rmk}

\begin{fact}  {\rm(\cite{BS})} $\frak{b}_\kappa$ is the least cardinality of any $F \subseteq {\cal C}_\kappa$ such that for any $A \in [\kappa]^\kappa$, there is $C \in F$ with $\vert A \setminus C \vert = \kappa$. 
\end{fact}

\begin{Obs}   \begin{enumerate}[\rm (i)]
\item $\frak{b}_\kappa$ is the least cardinality of any $F \subseteq {\cal C}_\kappa$ such that for any $D \in {\cal C}_\kappa$, there is $C \in F$ with $\vert D \setminus C \vert = \kappa$. 
\item There is a maximal descending $(I_\kappa, I_\kappa)$-tower of length ${\frak b}_\kappa$ consisting of closed unbounded subsets of $\kappa$.
\item Let $S \in NS^+_\kappa$ be such that $NS_\kappa \vert S$ is not $\sigma$-saturated, where $\sigma$ is an infinite cardinal less than or equal to $\frak{b}_\kappa$. Then $NS_\kappa \vert S$ is not $I_\kappa$-$\sigma$-saturated.
\end{enumerate}
\end{Obs}

{\bf Proof.} (i) : Let  $F \subseteq {\cal C}_\kappa$ with $0 < \vert F \vert\ < \frak{b}_\kappa$. By Fact 2.35, there must be $A \in [\kappa]^\kappa$ such that $\vert A \setminus C \vert < \kappa$ for all $C \in F$. Set $D = acc (A)$.
Then clearly, $D \in {\cal C}_\kappa$. Moreover, $\vert D \setminus C \vert < \kappa$ for all $C \in F$. 

(ii) : By Fact 2.35, we may find $D_i \in {\cal C}_\kappa$ for $i < {\frak b}_\kappa$ such that for any $A \in [\kappa]^\kappa$, there is $i < {\frak b}_\kappa$ with $\vert A \setminus D_i \vert = \kappa$. We inductively define $C_i \in {\cal C}_\kappa$ as follows. Set $C_0 = D_0$. Now suppose that $i > 0$ and $C_j$ has been constructed for each $j < i$. By (i), there is  $H \in {\cal C}_\kappa$ such that $\vert H \setminus C_j \vert < \kappa$ for all $j < i$. We let $C_i = H$ if $i$ is a limit ordinal, and $C_i = H \cap acc (C_{i - 1})$ otherwise. It is easy to see that $\langle C_i : i < {\frak b}_\kappa \rangle$ is a descending $(I_\kappa, I_\kappa)$-tower.
  
 (iii) : Select $A_\alpha \in (NS_\kappa \vert S)^+$ for $\alpha < \sigma$ such that $A_\beta \cap A_\alpha \in NS_\kappa \vert S$ whenever $\beta < \alpha < \sigma$. For $\beta < \alpha < \sigma$, pick $C_{\beta\alpha} \in {\cal C}_\kappa$ such that $(A_\beta \cap S) \cap (A_\alpha \cap S) \cap C_{\beta\alpha} = \emptyset$. For each $\alpha< \sigma$, there is by (i) $D_\alpha \in {\cal C}_\kappa$ such that $\vert D_\alpha \setminus C_{\beta\alpha} \vert < \kappa$ for all $\beta < \alpha$. Then clearly, $\vert (A_\beta \cap S \cap D_\beta) \cap (A_\alpha \cap S \cap D_\alpha) \vert < \kappa$ whenever $\beta < \alpha < \sigma$.
\hfill$\square$

\begin{Obs}  \begin{enumerate}[\rm (i)]
\item Let $J$ be a $\kappa$-complete ideal on $\kappa$. If there is a descending $(J, J)$-tower of length $\delta$, where $\delta \leq \kappa$, then there is a descending $(J, \{\emptyset\})$-tower of length $\delta$.
\item Let $J$ be a normal ideal on $\kappa$. If there is a descending $(J, J)$-tower of length $\delta$, where $\delta \leq \kappa^+$, then there is a descending $(J, I_\kappa)$-tower of length $\delta$.
\item Let $J = NS_\kappa \vert S$ for some $S \in NS^+_\kappa$. If there is a descending $(J, J)$-tower of length $\delta$, where $\delta \leq \frak{b}_\kappa$, then there is a descending $(J, I_\kappa)$-tower of length $\delta$.
\end{enumerate}
\end{Obs}

{\bf Proof.} We prove (iii) and leave the similar proofs of (i) and (ii) to the reader. Thus  suppose that $J = NS_\kappa \vert S$, where $S \in NS^+_\kappa$, and $\langle A_i : i < \delta \rangle$ is a descending $(J, J)$-tower of length $\delta$, where $0 < \delta \leq \kappa$. We recursively define $B_i \in J^+ \cap P (A_i \cap S)$ with $A_i \setminus B_i \in J$ for $i < \delta$ as follows. Put $B_0 = A_0$. Now suppose that $i > 0$, and $B_j$ has been constructed for every $j < i$. For each $j < i$, pick $C_{ji} \in {\cal C}_\kappa$ so that $(A_i \setminus B_j) \cap S \cap C_{ji} = \emptyset$. By Fact 2.35, there must be $D_i \in {\cal C}_\kappa$ such that $\vert D_i \setminus C_{ji} \vert < \kappa$ for all $j < i$. We set $B_i = A_i \cap S \cap D_i$. Notice that given $j < i$, $B_i \setminus B_j \subseteq D_i \setminus C_{ji}$, and therefore $\vert B_i \setminus B_j \vert < \kappa$. Furthermore, $B_j \setminus B_i \subseteq (A_j \setminus A_i) \setminus (A_j \setminus B_j)$, and consequently $B_j \setminus B_i \in J^+$.  
\hfill$\square$

\begin{Obs}  Let $J$ be a $\kappa$-complete ideal on $\kappa$, and $\sigma$ be a cardinal with $2 \leq \sigma \leq \kappa$. Then the following are equivalent :
 \begin{enumerate}[\rm (i)]
\item $J$ is not $\sigma$-saturated.
\item There is a descending $(J, \{\emptyset\})$-tower of length $\sigma$.
\item There is an ascending $(J, \{\emptyset\})$-tower of length $\sigma$.
\end{enumerate}
\end{Obs}

{\bf Proof.} (ii) $\rightarrow$ (i) and (iii) $\rightarrow$ (i) : By Observation 2.29.

(i) $\rightarrow$ (ii) and (iii) : Let $S_i \in J^+$ for $i < \sigma$ be such that $S_i \cap S_j \in J$ whenever $j < i < \sigma$. For $i < \sigma$, set $T_i = S_i \setminus (\bigcup_{j < i} S_j)$. Now we can define  an ascending $(J, \{\emptyset\})$-tower $\langle A_i : i < \sigma \rangle$ by $A_i = \bigcup_{j \leq i} T_j$, and a descending $(J, \{\emptyset\})$-tower $\langle B_i : i < \sigma \rangle$ by $B_i = \bigcup_{i \leq k < \sigma} T_k$.
\hfill$\square$

\begin{Obs}   \begin{enumerate}[\rm (i)]
\item Let $J$ be a $\kappa$-complete ideal on $\kappa$, $\delta$ be a nonzero ordinal, and $\langle A_i : i < \delta \rangle$ be a maximal ascending $(J, Y)$-tower, where $Y\subseteq P(\kappa)$. Then $\delta$ is a successor ordinal.
\item Let $J$ be a $\kappa$-complete, nowhere prime ideal on $\kappa$, $\delta$ be a nonzero ordinal, and $\langle A_i : i < \delta \rangle$ be a maximal descending $(J, Y)$-tower, where $Y$ is a subset of $P(\kappa)$ closed under subsets. Then $\delta$ is not a successor ordinal.
\end{enumerate}
\end{Obs}

{\bf Proof.} (i) : Suppose otherwise. Put $A_\delta = \bigcup_{i < \delta} A_i$. Then $\langle A_j : j \leq \delta \rangle$ is an ascending $(J, Y)$-tower. Contradiction.

(ii) : Suppose otherwise, and let $\delta = \xi + 1$. Then $A_\xi$ can be written as the disjoint union of two members of $J^+$, say $B_0$ and $B_1$. Put $A_\delta = B_0$. Then $\langle A_j : j \leq \delta \rangle$ is a descending $(J, Y)$-tower. Contradiction.
\hfill$\square$

\begin{Obs}  Let $J$ be a $\kappa$-complete ideal on $\kappa$, and $\sigma$ be a cardinal with $2 \leq \sigma \leq \kappa$ such that $J$ is not $\sigma$-saturated. Then the following hold :
 \begin{enumerate}[\rm (i)]
\item Suppose $\sigma \geq \omega$. Then there is a maximal descending $(J, \{\emptyset\})$-tower of length $\sigma$.
\item There is a maximal ascending $(J, \{\emptyset\})$-tower of length $\delta$, where $\delta$ equals $\sigma$ if $\sigma < \omega$, and $\sigma + 1$ otherwise.
\end{enumerate}
\end{Obs}

{\bf Proof.} Use (the proof of) Observations 2.38 and 2.39.
\hfill$\square$

\begin{Obs}   Let $J$ be a normal ideal on $\kappa$ that is not $\kappa$-saturated. Then there is a maximal descending $(J, J)$-tower of length $\kappa$.
\end{Obs}

{\bf Proof.}  We use an argument that can be found in \cite{Tall}. Pick a partition $\langle A_\alpha : \alpha < \kappa \rangle$ of $\kappa$ into members of $J^+$. Set $B_0 = \bigcup_{0 > \alpha < \kappa} (A_\alpha \cap (\alpha + 1))$, and $B_\alpha = A_\alpha \setminus (\alpha + 1)$ for $0 < \alpha < \kappa$. Put $S_i = \bigcup_{\alpha \geq i} B_\alpha$ for each $i < \kappa$. Now let $S \subseteq \kappa$ be such that $S \setminus S_i \in J$. for all $i < \kappa$. For each $\alpha < \kappa$, we may find $C_\alpha \in J^\ast$ such that $S \cap B_\alpha \cap C_\alpha = \emptyset$. Then 

\centerline{$S \cap \bigtriangleup_{\alpha < \kappa} C_\alpha= \bigcup_{\alpha < \kappa} (S \cap (B_\alpha \cap \bigtriangleup_{\alpha < \kappa} C_\alpha)) \subseteq 1$,}

and therefore $S \in J$. Thus $\langle S_i : i < \kappa \rangle$ is a maximal descending $(J, J)$-tower.
\hfill$\square$

\medskip

The following is due to Moti Gitik \cite{Gitik}.

\medskip
 
\begin{Th}   Let $J$ be a normal ideal on $\kappa$ that is not $\kappa^+$-saturated. Then there is a descending $(J, J)$-tower of length $\kappa^+$.
\end{Th}

{\bf Proof.}  Pick $A_\alpha \in J^+$ for $\alpha < \kappa^+$ such that $A_\beta \cap A_\alpha \in J$ whenever $\beta < \alpha < \kappa^+$. For $\kappa \leq \alpha < \kappa^+$, select a bijection $j_\alpha : \kappa \rightarrow \alpha$ and put $B_\alpha = \bigtriangleup_{i < \kappa} (\kappa \setminus A_{j_\alpha (i)})$. Note that $\vert A_\beta \cap B_\alpha \vert < \kappa$ for all $\beta < \alpha$.

\medskip

{\bf Claim 1.} Let $\kappa \leq \beta < \kappa^+$. Then $A_\beta \setminus B_\beta \in J$. 

\medskip

{\bf Proof of Claim 1.}  Suppose otherwise. Define $f : A_\beta \setminus B_\beta \rightarrow \kappa$ by $f (\xi) =$ the least $i < \xi$ such that $\xi \in A_{j_\beta (i)}$. There must be $H \in J^+ \cap P (A_\beta \setminus B_\beta)$ such that $f$ is constant on $H$. This contradiction completes the proof of the claim.

\medskip

{\bf Claim 2.} Let $\kappa \leq \beta < \alpha < \kappa^+$. Then $B_\beta \setminus B_\alpha \in J^+$. 

\medskip

{\bf Proof of Claim 2.}  Clearly, $A_\beta \setminus (B_\beta \setminus B_\alpha)$ is a subset of $(A_\beta \setminus B_\beta) \cup (A\beta \cap B_\alpha)$ which by Claim 1 lies in $J$. Hence $B_\beta \setminus B_\alpha \in J^+$, which completes the proof of the claim.

\medskip

{\bf Claim 3.} Let $\kappa \leq \beta < \alpha < \kappa^+$. Then $B_\alpha \setminus B_\beta \in J$. 

\medskip

{\bf Proof of Claim 3.}  Suppose otherwise. Define $g : B_\alpha \setminus B_\beta \rightarrow \kappa$ by $g (\xi) =$ the least $i < \xi$ such that $\xi \in A_{j_\beta (i)}$. We may find $G \in J^+ \cap P (B_\alpha \setminus B_\beta)$ and $i < \kappa$ such that $g$ takes the constant value $i$ on $G$. Let $k < \kappa$ be such that $j_\beta (i) = j_\alpha (k)$. Then $\xi \leq k$ for all $\xi \in G$. This contradiction completes the proof of the claim and that of the proposition.
\hfill$\square$

\begin{Rmk} Suppose that in the proof above, the family $\{A_\alpha : \alpha < \kappa^+\}$ has the additional property that for any $K \in J^+$, there is $\alpha$ with $K \cap A_\alpha \in J^+$. Then our $(J, J)$-tower $\langle B_\beta : \kappa \leq \beta < \kappa^+ \rangle$ is maximal. To see this, recall that $\vert B_\beta \cap A_\alpha \vert < \kappa$ whenever $\alpha < \beta$ and $\kappa \leq \beta < \kappa^+$. It follows that if $W \subseteq \kappa$ is such that $W \setminus  B_\beta \in J$ whenever $\kappa \leq \beta < \kappa^+$, then $W \in J$.
\end{Rmk}

\begin{Def}  $Depth ({}^\kappa \kappa)$  denotes the least ordinal $\eta$ such that there is no increasing sequence $\langle f_i : i < \eta \rangle$ in $({}^\kappa \kappa, <^\ast)$.
\end{Def}

\begin{Def}  $Depth ({\cal C}_\kappa)$  denotes the least ordinal $\delta$ such that there is no sequence $\langle C_i : i < \delta \rangle$ such that
\begin{itemize}
\item $C_i \in {\cal C}_\kappa$ ;
\item $C_{i + 1} \subseteq acc (C_i)$ ;
\item $\vert C_i \setminus C_j \vert < \kappa$ for all $j < i$.
\end{itemize}
\end{Def}

\begin{Obs}   $Depth ({}^\kappa \kappa)$ (respectively,  $Depth ({\cal C}_\kappa)$) is not the successor of a successor ordinal.
\end{Obs}

\begin{fact}  {\rm(\cite{She589})} $\frak{b}_\kappa < Depth ({\cal C}_\kappa) \leq Depth ({}^\kappa \kappa) \leq Depth ({\cal C}_\kappa) + 1$.
\end{fact}

{\bf Proof.} For the first inequality see the proof of Observation 2.36 (ii). To establish the second one, let $\langle C_i : i < \delta \rangle$ be such that
\begin{itemize}
\item $C_i \in {\cal C}_\kappa$ ;
\item $C_{i + 1} \subseteq acc (C_i)$ ;
\item $\vert C_i \setminus C_j \vert < \kappa$ for all $j < i$.
\end{itemize}
For $S \subseteq \kappa$, let $e_S : o.t. (S) \rightarrow S$ be the increasing enumeration of $S$. For $i < \delta$, define $f_i \in {}^\kappa \kappa$ by $f_i (\alpha) = e_{C_i} (\alpha + 1)$. Now fix $i < \delta$. Then for any $\beta < \kappa$, $e_{C_i} (\beta) \leq e_{acc (C_i)} (\beta)$. Hence for each $\alpha < \kappa$,

\centerline{$f_i (\alpha) < e_{C_i} (\alpha + \omega) \leq e_{acc (C_i)} (\alpha + 1) \leq f_{i + 1} (\alpha)$.}

Given $i + 1 < j < \delta$, we may find $\xi < \eta < \kappa$ such that $C_j \setminus C_{i + 1} \subseteq \xi$ and $o.t. (C_{i + 1} \cap \eta) = \eta = o.t. (C_j \cap \eta)$. Then clearly,

\centerline{$f_i (\alpha) < f_{i + 1} (\alpha) \leq f_j (\alpha)$}

whenever $\eta \leq \alpha < \kappa$. 

For the last inequality, given an increasing sequence $\langle f_i : i \leq \delta \rangle$ in $({}^\kappa \kappa, <^\ast)$, let

$D$ denote  the set of all $\alpha \in acc (\kappa)$ such that 
\begin{itemize}
\item   $f_\delta (\beta) < \alpha$ for all $\beta < \alpha$ ;
\item $\omega^\alpha = \alpha$.
\end{itemize}
Now for each $i < \delta$, set

\centerline{$C_i = \{ \alpha + \omega^{f_i (\alpha)} \cdot \gamma : \alpha \in D$ and $\gamma < f_i (\alpha) < f_\delta (\alpha)\}$.}

It is not difficult to see that

\begin{itemize}
\item $C_i \in {\cal C}_\kappa$ ;
\item $C_{i + 1} \subseteq acc (C_i)$ ;
\item $\vert C_i \setminus C_j \vert < \kappa$ for all $j < i$.
\end{itemize}

\hfill$\square$

\begin{Def}  For $f \in \kappa^\kappa$, let $M_f = \{(\alpha, \beta) \in \kappa \times \kappa : \beta \geq f (\alpha)\}$ and   $m_f = \{(\alpha, \beta) \in \kappa \times \kappa : \beta \leq f (\alpha)\}$. 
\end{Def}

\begin{Rmk} Since $f \subseteq m_f \cap M_f$, we have $\{m_f, M_f\} \subset [\kappa \times \kappa]^\kappa$.
\end{Rmk}

\begin{Pro}  Let $f, g \in \kappa^\kappa$ be such that $f <^\ast g$. Then the following hold : 
\begin{enumerate}[\rm (i)]
\item $\vert M_g \setminus M_f \vert < \kappa$, and moreover $\vert M_f \setminus M_g \vert = \kappa$.
\item $\vert m_f \setminus m_g \vert < \kappa$, and moreover $\vert m_g \setminus m_f \vert = \kappa$.
\end{enumerate}
\end{Pro}

{\bf Proof.} Let $\gamma < \kappa$ be such that $f (\alpha) < g (\alpha)$ for all $\alpha \geq \gamma$. Then the following is readily checked :
\begin{itemize}
\item $M_g \setminus M_f \subseteq \bigcup_{\alpha < \gamma} \{ (\alpha, \beta) : \beta < f (\alpha)\}$.
\item $\{(\alpha, f (\alpha)) : \alpha \geq \gamma \} \subseteq M_f \setminus M_g$.
\item $m_f \setminus m_g \subseteq \bigcup_{\alpha < \gamma} \{ (\alpha, \beta) : g (\alpha) < \beta \leq f (\alpha)\}$.
\item $\{(\alpha, g (\alpha)) : \alpha \geq \gamma \} \subseteq m_g \setminus m_f$.
\end{itemize} 
\hfill$\square$

\begin{Cor}   \begin{enumerate}[\rm (i)]
\item   {\rm(\cite{She589})} $Depth ({}^\kappa \kappa) \leq Depth ([\kappa]^\kappa, \nearrow)$.
\item  $Depth ({}^\kappa \kappa) \leq Depth ([\kappa]^\kappa, \searrow)$.
\end{enumerate}
\end{Cor}

\begin{Rmk} By Theorem 2.31, Fact 2.47 and Corollary 2.51, $\clubsuit^{\rm ev}_{\omega_1} [NS_{\omega_1}]$ implies the existence of a descending (respectively, ascending) $(NS_{\omega_1}, I_ {\omega_1})$-tower of length ${\frak b}_{\omega_1}$. Let us mention in this connection that by a result of Baumgartner and Tall \cite{Tall}, $\diamondsuit_{\omega_1} [NS_{\omega_1}]$ (and hence $\clubsuit^{\rm ev}_{\omega_1} [NS_{\omega_1}]$) + $2^{\aleph_1} > \aleph_2$ + $P_S$ is consistent, where $P_S$ asserts the following :

For any family $F$ of less than $2^{\aleph_1}$ many stationary subsets of $\omega_1$ with the property that $\bigtriangleup_{\alpha < \omega_1} f (\alpha) \in NS^+_{\omega_1}$ for all $f : \omega_1 \rightarrow F$, there is a stationary subset $T$ of $\omega_1$ such that $\vert T \setminus S \vert < \kappa$ for every $S \in F$. 

Notice that it is immediate from Fact 2.35 that $P_S$ implies ${\frak b}_{\omega_1} = 2^{\aleph_1}$.
\end{Rmk}


\begin{Rmk} Suppose that the GCH holds in $V$, and let $\tau$ and $\sigma$ be two regular cardinals with $\kappa^+ \leq \tau < \sigma$. Put ${\mathbb Q} = (\{0\} \times \tau) \cup (\{1\} \times \sigma)$. For $q, r \in{\mathbb Q} $, let $q < r$ just in case either $q = (0, \alpha)$ and $r = (0, \beta)$, where $\alpha < \beta < \tau$, or $q = (1, \gamma)$ and $r = (1, \delta)$, where $\gamma < \delta < \sigma$. Notice that $\tau =$ the least size of any unbounded subset of ${\mathbb Q}$. Furthermore ${\mathbb Q}$ is well-founded. Hence by a result of Cummings and Shelah \cite{CS}, there is a $\kappa$-closed, $\kappa^+$-cc notion of forcing ${\mathbb P}$ such that in $V^{\mathbb{P}}$,
\begin{itemize}
\item ${\frak b}_\kappa = \tau$.
\item There are $f_q \in {}^\kappa \kappa$ for $q \in {\mathbb Q}$ such that  
\begin{enumerate}[\rm (a)]
\item for any $g \in {}^\kappa \kappa$, there is $q \in {\mathbb Q}$ with $g <^\ast f_q$ ;
\item for $q, r \in {\mathbb Q}$, $q < r$ if and only if $f_q <^\ast f_r$.
\end{enumerate}
\end{itemize}
Thus in  $V^{\mathbb{P}}$, $Depth ({}^\kappa \kappa) > \sigma$.
\end{Rmk}

\begin{Rmk} Gitik observes that there is a natural forcing (let $P$ be the set of all  $(c, F)$ such that $c$ is a closed subset of $\kappa$ of size less than $\kappa$, and $F \in P_\kappa ({\cal C}_\kappa)$, with the obvious ordering)  that adds $C \in {\cal C}_\kappa$ such that $\vert C \setminus D \vert < \kappa$ for every $D$ in $({\cal C}_\kappa)^V$. It can be iterated to any length, which tends to indicate that there is no nontrivial upper bound for $Depth ({\cal C}_\kappa)$.
\end{Rmk}

\bigskip

\subsection{Interdependent depths}

\medskip

Let us next discuss the following result of Shelah. 

\medskip

\begin{fact}  {\rm (\cite{She589})} Suppose that
\begin{itemize}
\item $\rho$ is a singular cardinal of uncountable cofinality such that $\sigma^{\cf (\rho)} < \rho$ for every cardinal $\sigma < \rho$ ;
\item $\langle \rho_i : i < \cf (\rho) \rangle$ is an increasing, continuous sequence of infinite cardinals with supremum $\rho$ ;
 \item $f$ is a function from $\cf (\rho)$ to the set of all regular infinite cardinals below $\rho$ such that $f (i) < Depth ({\cal C}_{\rho_i^+})$ for all $i < \cf (\rho)$ ;
 \item $I$ is a normal ideal on $\cf (\rho)$ ;
 \item $\pi = tcf (\prod f / I)$.
\end{itemize}
Then $\pi < Depth ({\cal C}_{\rho^+})$.
 \end{fact}
\bigskip

\begin{Rmk} \cite{She589} contains more results of the same type.
\end{Rmk}

\medskip

Let us consider a concrete situation where Fact 2.55 can be applied. In \cite{Mer}  Merimovich constructs from large large cardinals a number of models where GCH massively fails. To be specific let us choose a model $V^{\mathbb{P}}$ in which there are an inaccessible cardinal $\theta$ and $C$ in ${\cal C}_\theta$ consisting of infinite cardinals such that for any infinite cardinal $\tau < \theta$, $2^\tau$ equals $\sigma^{+ 3}$ if there is $\sigma \in acc (C)$ such that $\sigma \leq \tau < \sigma^{+ 3}$, and $\tau^+$ otherwise. As pointed out by Gitik to the author, it can be arranged that in $V^{\mathbb{P}}$ , $\frak{b}_{\sigma^+} = 2^\sigma$ for every $\sigma \in acc (C)$. Now working in  $V^{\mathbb{P}}$, let $\rho \in acc (C)$ be a singular cardinal of uncountable cofinality, and let $\langle \rho_i : i < \cf (\rho) \rangle$ be an increasing, continuous sequence of singular cardinals in $acc (C)$ with supremum $\rho$. Then $pp (\rho) = 2^\rho = \rho^{+ 3}$. Hence by Lemma 9.2.9 in \cite{HSW}, there must be $D \in {\cal C}_{\cf (\rho)}$ and a function $f$ from $D$ to the set of all regular infinite cardinals below $\rho$ such that 
\begin{itemize}
 \item $ tcf (\prod f / J) = \rho^{+ 3}$, where $J$ denotes the noncofinal ideal on $C$ ;
\item for any $i \in D$,

\centerline{$\rho_i < f (i) \leq pp_{\cf (\rho)} (\rho_i) \leq 2^{\rho_i} = \rho_i^{+ 3}.$}

 \end{itemize}
 Let $I$ be the nonstationary ideal on $C$. Then by Lemma 3.17 of \cite{HSW}, $ tcf (\prod f / I) =  tcf (\prod f / J) = \rho^{+ 3}$. Hence by Fact 2.55, $Depth ({\cal C}_{\rho^+}) > 2^{(\rho^+)}$. 
 
 The interpretation is that $Depth ({\cal C}_{\rho^+})$ depends on the $Depth ({\cal C}_{\rho_i^+})$'s . But what about $\frak{b}_{\rho^+}$ in such a situation ? Does it also depend on the $\frak{b}_{\rho_i^+}$'s ?
 
 \bigskip

\section{J'ENLEVE TOUT}

\medskip

\subsection{Fromage ou dessert}

\medskip

Let us return to our starting point, when $\kappa = \nu^+$, $\theta$ is a regular cardinal less than $\nu$, and $J$ is a $\kappa$-complete ideal on $\kappa$ extending $NS_\kappa \vert E^\kappa_\theta$. We just saw that if  $\clubsuit^{\rm ev}_\kappa [J]$ holds, then there is a descending (respectively, ascending) $(J, I_ \kappa)$-tower of length ${\frak b}_\kappa$. But what if $\clubsuit^{\rm ev}_\kappa [J]$ fails ? To address this problem, we could weaken our club principle in a number of ways. We could for instance allow several guesses instead of just one. 

\medskip

\begin{Def}  Given a $\kappa$-complete ideal $J$ on $\kappa$ and a cardinal $\tau$ with $1 \leq \tau < \kappa$, $\clubsuit_\kappa^{- / \tau} [J]$ asserts the existence of $B^{i}_\delta \subseteq \delta$ with $\sup B^{i}_\delta = \delta$  for $\delta \in acc (\kappa)$ and $i < \tau$ such that 

\centerline{$\{ \delta \in acc (\kappa) : \exists i < \tau (B^{i}_\delta \subseteq W)\} \in J^+$}

 for any $W \in [\kappa]^\kappa$.

$\clubsuit_\kappa^- [J]$ asserts the existence of $B^{i}_\delta \subseteq \delta$ with $\sup B^{i}_\delta = \delta$  for $\delta \in acc (\kappa)$ and $i < \delta$ such that 

\centerline{$\{ \delta \in acc (\kappa) : \exists i < \delta (B^{i}_\delta \subseteq W)\} \in J^+$}

 for any $W \in [\kappa]^\kappa$.
\end{Def}

\medskip

Note that $\clubsuit^{\rm ev}_\kappa [J] \Rightarrow \clubsuit_\kappa^-  [J]$. Furthermore, if $E^\kappa_\tau \in J^\ast$, then $\clubsuit^{\rm ev}_\kappa [J] \Rightarrow \clubsuit_\kappa^{- / \tau} [J]$. 

\medskip

\begin{Obs} \begin{enumerate}[(i)] 
\item Let $\langle B^{i}_\delta : \delta \in acc (\kappa) \rangle$ witness that  $\clubsuit_\kappa^{- / \tau} [J]$ holds, where $J$ is a $\kappa$-complete ideal on $\kappa$, and $\tau$ a cardinal with $1 \leq \tau < \kappa$. Then for any $W \in [\kappa]^\kappa$, there is $i < \tau$ such that  $\{ \delta \in acc (\kappa) : B^{i}_\delta \subseteq W\} \in J^+$.
\item Let $\langle B^{i}_\delta : \delta \in acc (\kappa) \rangle$ witness that  $\clubsuit_\kappa^- [J]$ holds, where $J$ is a normal ideal on $\kappa$. Then for any $W \in [\kappa]^\kappa$, there is $i < \kappa$ such that  $\{ \delta \in acc (\kappa) \setminus (i + 1) : B^{i}_\delta \subseteq W\} \in J^+$.
\end{enumerate} 
\end{Obs}

\begin{Obs} \begin{enumerate}[(i)] 
\item Suppose that  $\clubsuit_\kappa^{- / \tau} [J]$ holds, where $J$ is a $\kappa$-complete ideal on $\kappa$, and $\tau$ a cardinal with $1 \leq \tau < \kappa$, and let $\rho$ be a cardinal such that $\cf (\rho) > \tau$ and $I_\kappa$ is not $\rho$-saturated. Then $J$ is not $I_\kappa$-$\rho$-saturated.
\item Suppose that  $\clubsuit_\kappa^- [J]$ holds, where $J$ is a normal ideal on $\kappa$, and let $\rho$ be a cardinal such that $\cf (\rho) > \kappa$ and $I_\kappa$ is not $\rho$-saturated. Then $J$ is not $I_\kappa$-$\rho$-saturated.
\end{enumerate} 
\end{Obs}

{\bf Proof.} By Observation 3.2 and the proof of Observation 2.27.
\hfill$\square$

\begin{Obs} Given a $\kappa$-complete ideal $J$ on $\kappa$ extending $NS_\kappa$, the following are equivalent :
 \begin{enumerate}[\rm (i)]
 \item $\diamondsuit_\kappa^- [J]$ holds.
\item $\clubsuit_\kappa^- [J]$ holds and $2^{< \kappa} = \kappa$.
\end{enumerate}   
 \end{Obs}

{\bf Proof.}  (i) $\rightarrow$ (ii) : Use Fact 2.16 (ii).

(ii)$ \rightarrow$ (i) : Let  $B^{i}_\delta \subseteq \delta$ with $\sup B^{i}_\delta = \delta$  for $\delta \in acc (\kappa)$ and $i < \delta$ be such that $\{ \delta \in acc (\kappa) : \exists i < \delta (B^{i}_\delta \subseteq W)\} \in J^+$ for any $W \in [\kappa]^\kappa$. For $\eta < \kappa$, define $\chi_\eta : P (\eta) \rightarrow {}^\eta 2$ by : $\chi (\eta) (a) (\xi) = 1$ if and only if $\xi \in a$. Select a bijection $F : \bigcup_{\eta < \kappa} {}^\eta 2 \rightarrow \kappa$. For $\delta \in acc (\kappa)$ and $i < \delta$, put

\centerline{$s^{i}_\delta = \bigcup_{\eta < \kappa} (\bigcup\{a \subseteq \eta : F (\chi_\eta (a)) \in B^{i}_\delta\})$.}

Given $A \subseteq \kappa$, let $D$ be the set of all $\delta \in acc (\kappa)$ such that $\sup \{ F (\chi_\zeta (A \cap \zeta)) : \zeta \leq \alpha\} < \delta$ for all $\alpha < \delta$. Note that $D$ belongs to $NS_\kappa^\ast$ and hence to $J^\ast$. Now suppose that $\delta \in D$ and $i < \delta$ are such that $B^{i}_\delta \subseteq \{ F (\chi_\eta (A \cap \eta)) : \eta < \kappa \}$.

\medskip

{\bf Claim.}  $s^{i}_\delta \cap \delta = A \cap \delta$.

\medskip

{\bf Proof of the claim.}  

$\subseteq$ : Let $\alpha \in s^{i}_\delta \cap \delta$. We may find $\eta < \kappa$ and $a \subseteq \eta$ such that $\alpha \in a$ and $F (\chi_\eta (a)) \in B^{i}_\delta$. Then clearly, $F (\chi_\eta (a)) = F (\chi_\eta (A \cap \eta))$, and therefore $a = A \cap \eta$. Hence, $\alpha \in A \cap \delta$.

$\supseteq$ : Let $\alpha \in A \cap \delta$. There must be $\gamma \in B^{i}_\delta$ such that $\gamma > \sup \{ F (\chi_\zeta (A \cap \zeta)) : \zeta \leq \alpha\}$. Let $\eta < \kappa$ be such that $\gamma = F (\chi_\eta (A \cap \eta))$. Since $\eta > \alpha$, we have that $\alpha \in A \cap \eta$. It follows that $\alpha \in s^{i}_\delta$, which completes the proof of the claim and that of the observation. 
\hfill$\square$

\medskip

There is another way to weaken  $\clubsuit^{\rm ev}_\kappa [J]$. Instead of guessing eventually,  we could content ourselves with guessing cofinally. But then we need an extra condition on our guess $B_\delta$, otherwise we would achieve success too easily with $B_\delta = \delta$. 

\medskip

\begin{Def}  Given a $\kappa$-complete ideal $J$ on $\kappa$ and a cardinal $\sigma < \kappa$, $\clubsuit^{\rm cof / \sigma}_\kappa [J]$ asserts the existence of $B_\delta \in P_\sigma (\delta)$  for $\delta < \kappa$ such that $\{ \delta : \sup (W \cap B_\delta) = \delta\} \in J^+$ for any $W \in [\kappa]^\kappa$.

$\clubsuit^{\rm cof}_\kappa [J]$ asserts the existence of $B_\delta \in P_{\vert \delta \vert} (\delta)$  for $0 < \delta < \kappa$ such that $\{ \delta : \sup (W \cap B_\delta) = \delta\} \in J^+$ for any $W \in [\kappa]^\kappa$.\end{Def}

\medskip

Note that if $E^\kappa_\theta \in J^\ast$, where $\theta^+ < \kappa$, then $\clubsuit^{\rm ev}_\kappa [J] \Rightarrow \clubsuit^{\rm cof / \theta^+}_\kappa [J]$. 

\medskip

We finally settle for a doubly weaker principle.

\medskip

\begin{Def}  Let $\sigma$ and $\tau$ be two cardinals with $\sigma < \kappa$ and $1 \leq \tau < \kappa$, and $J$ be a $\kappa$-complete ideal on $\kappa$. The principle $\clubsuit_\kappa^{{\rm cof / \sigma}, - / \tau} [J]$ asserts the existence of $B^{i}_\delta \in P_\sigma (\delta)$  for $\delta < \kappa$ and $i < \tau$ such that for any $W \in [\kappa]^\kappa$,

\centerline{$\{ \delta < \kappa : \exists i < \tau (\sup (W \cap B^{i}_\delta) = \delta)\} \in J^+$.}

$\clubsuit_\kappa^{{\rm cof / \sigma}, -} [J]$ (respectively, $\clubsuit_\kappa^{{\rm cof}, -} [J]$) asserts the existence of $B^{i}_\delta$ in $P_\sigma (\delta)$ (respectively, $P_{\vert \delta \vert} (\delta)$)  for $i < \delta < \kappa$ such that for any $W \in [\kappa]^\kappa$,

\centerline{$\{ \delta < \kappa : \exists i < \delta (\sup (W \cap B^{i}_\delta) = \delta)\} \in J^+$.}
\end{Def}

\medskip

It is easy to see that if $\clubsuit_\kappa^{{\rm cof}, -} [J]$ holds, then $NS^\ast_\kappa \subseteq J^+$. Notice that if $\kappa = \nu^+$, then $\clubsuit_\kappa^{{\rm cof / \sigma}, - / \nu} [J]$ (respectively, $\clubsuit_\kappa^{{\rm cof / \nu}, -} [J]$, $\clubsuit_\kappa^{- / \nu} [J]$) and $\clubsuit_\kappa^{{\rm cof / \sigma}, -} [J]$ (respectively, $\clubsuit_\kappa^{{\rm cof}, -} [J]$, $\clubsuit_\kappa^- [J]$) are equivalent. 

\medskip

\begin{Obs}  Suppose that $\kappa = \nu^+$ and $\clubsuit_\kappa^{{\rm cof}, -} [J]$ holds. Then the following hold :
\begin{enumerate}[\rm (i)]
\item $\nu > \omega$.
\item  $E^\kappa_{< \nu} \in J^+$, and moreover $\clubsuit_\kappa^{{\rm cof}, -} [J \vert E^\kappa_{< \nu}]$ holds.
\item ${\cal C}_\kappa \subseteq J^+$.
\end{enumerate}
\end{Obs}

\medskip

Assuming $\kappa$ is a successor cardinal, $\clubsuit_\kappa^{{\rm cof}, -} [NS_\kappa \vert S]$ is denoted by $\clubsuit^-_S$ in \cite {Rin10} where, building on previous work by D\v zamonja and Shelah \cite{DZ1}, Rinot proves that  $\clubsuit^-_S$ implies that $NS_\kappa \vert S$ is not $\kappa^+$-saturated. Our presentation will closely follow his. We start with the following technical lemma.

\medskip

\begin{Lem}  Suppose that $\kappa = \nu^+$, and $J$ is a $\kappa$-complete ideal on $\kappa$ such that $\clubsuit_\kappa^{{\rm cof}, -} [J]$ holds. Then there exist $A^{i}_\delta \in P_\nu (\kappa \times \kappa)$ for $\delta < \kappa$ and $i < \nu$ such that for any $F : \kappa \rightarrow \kappa$,

\centerline{$\{ \delta < \kappa : \exists i < \nu \exists Z \subseteq \delta (\sup Z = \delta$ and $F \vert Z \subseteq A^{i}_\delta)\}$}

lies in $J^+$. 
\end{Lem}

{\bf Proof.} We follow the proof of Lemma 1.5 in \cite{Rin10}. Let ${\cal B}_\delta \in P_\kappa (P_\nu (\delta))$ for $\delta < \kappa$ witness that $\clubsuit_\kappa^{{\rm cof}, -} [J]$ holds. Set $\chi = \cf (\nu)$. Let $\langle \xi_j : j < \chi \rangle$ be an increasing sequence of infinite ordinals with supremum $\nu$, and for each $\gamma \in E^\kappa_{\not= \chi}$, select $G^j_\gamma \subseteq \gamma \times \gamma$ for $j < \chi$ so that
\begin{itemize}
\item $\vert G^j_\gamma \vert \leq \vert \xi_j \vert$ ;
\item $G^k_\gamma \subseteq G^j_\gamma$ for all $k < j$ ;
\item $ \bigcup_{j < \chi} G^j_\gamma = \gamma \times \gamma$.
\end{itemize}
For $\delta < \kappa$, let $\langle A^{i}_\delta : i < \nu \rangle$ be an enumeration of the set

\centerline{$\{A :  \exists B \in {\cal B}_\delta  \exists j < \chi (A = \bigcup \{G^j_\gamma : \gamma \in B \cap E^\kappa_{\not= \chi}\})\}$.}

Now fix $F : \kappa \rightarrow \kappa$. Put $X = \{\gamma \in E^\kappa_{\not= \chi} : \forall \alpha < \gamma (F (\alpha) < \gamma) \}$.

\medskip

{\bf Claim 1.}  Let $\gamma \in X$. Then $F \vert Y \subseteq G^j_\gamma$ for some cofinal subset $Y$ of $\gamma$ and some $j < \chi$.

\medskip

{\bf Proof of Claim 1.}  Pick a cofinal subset $e$ of $\gamma$ of order-type $\cf (\gamma)$. Define $h : e \rightarrow \chi$ by $h (\alpha) = \min\{j < \chi : (\alpha, F (\alpha)) \in G^j_\gamma\}$. Then we may find $Y \subseteq e$ with
$\vert Y \vert = \cf (\gamma)$, and $j < \chi$ such that $h (\alpha) \leq j$ for all $\alpha \in Y$. Clearly, $Y$ and $j$ are as desired, which completes the proof of the claim. 

\medskip

Using Claim 1, define $g : X \rightarrow \chi$ by $g (\gamma) =$ the least $j$ such that $F \vert Y \subseteq G^j_\gamma$ for some cofinal subset $Y$ of $\gamma$. Then we may find $j < \chi$ and a size $\kappa$ subset $W$ of $X$ such that $g$ takes the constant value $j$ on $W$. Then $S \in J^+$, where $S = \{ \delta < \kappa : \exists B \in {\cal B}_\delta (\sup (W \cap B) = \delta)\}$.

\medskip

{\bf Claim 2.}  Let $\delta \in S$, and let $B \in {\cal B}_\delta$ such that $\sup (W \cap B) = \delta$. Then there exists a cofinal subset $Z$ of $\delta$ such that $F \vert Z \subseteq  \bigcup \{G^j_\gamma : \gamma \in B \cap E^\kappa_{\not= \chi}\}$.

\medskip

{\bf Proof of Claim 2.}  For each $\gamma \in Z \cap B$, select a cofinal subset $Y_\gamma$ of $\gamma$ with $F \vert Y_\gamma \subseteq G^j_\gamma$. Then $W = \bigcup \{Y_\gamma : \gamma \in Z \cap B \}$ is as desired, 
 which completes the proof of the claim and that of the lemma. 

\hfill$\square$

\begin{Th}  Suppose that 
\begin{itemize}
\item $\kappa = \nu^+$  ;
\item $J$ is a $\kappa$-complete ideal on $\kappa$ such that $\clubsuit_\kappa^{{\rm cof}, -} [J]$ holds ; 
\item $\tau$ is a regular cardinal less than $Depth ({}^\kappa \kappa)$.
\end{itemize}
Then there exists either an ascending $(J, I_\kappa)$-tower of length $\tau$, or a descending $(J, J)$-tower of length $\tau$.
\end{Th}

{\bf Proof.} The proof is a modification of that of Theorem 1.10 in \cite{Rin10}. Select $f_\alpha \in {}^\kappa \kappa$ for $\alpha < \tau$ such that $f_\alpha <^\ast f_\beta$ whenever $\alpha < \beta < \tau$. Let $A^{i}_\delta \in P_\nu (\kappa \times \kappa)$ for $\delta < \kappa$ and $i < \nu$ be as in the statement of Lemma 3.8. For each $\alpha < \tau$, pick $i_\alpha < \nu$ such that

\centerline{$S_\alpha = \{ \delta \in E^\kappa_{< \nu} :  \exists W \subseteq \delta (\sup W = \delta$ and $f_\alpha \vert W \subseteq A^{i_\alpha}_\delta)\}$}

lies in $J^+$. By thinning out our sequence of functions, we may assume that there is $i < \nu$ such that $i_\alpha = i$ for all $\alpha < \tau$. 

For $\delta \in E^\kappa_{< \nu}$, put

\centerline{$D_\delta = \{j < \delta : \exists r ((j, r) \in A^{i}_\delta)\}$}

and

\centerline{$R_\delta = \{r < \delta : \exists j ((j, r) \in A^{i}_\delta)\}$.}

For $\delta \in E^\kappa_{< \nu}$ and $\alpha < \tau$, define $f_{\alpha\delta} : D_\delta \rightarrow R_\delta$ by $f_{\alpha\delta} (j) = \min ((R_\delta \cup \{\kappa\}) \setminus f_\alpha (j))$. Finally, for $\alpha < \beta < \tau$, let

\centerline{$S_{\alpha\beta} = \{\delta \in E^\kappa_{< \nu} :  \sup \{j \in D_\delta : f_{\alpha\delta} (j) < f_{\beta\delta} (j)\} = \delta\}.$}

\medskip

{\bf Claim 1.}  Let $\alpha < \beta < \tau$. Then $\vert S_\alpha \setminus S_{\alpha\beta} \vert < \kappa$ (and hence $S_{\alpha\beta} \in J^+$).

\medskip

{\bf Proof of Claim 1.}  Suppose not. Put $m = \sup \{j < \kappa : f_\alpha (j) \geq f_\beta (j)\}$, and pick $\delta \in S_\alpha \setminus S_{\alpha\beta}$ with $\delta > m$. Note that for any $j \in D_\delta$ with $j > m$, we have $f_\alpha (j) < f_\beta (j)$ and hence $f_{\alpha\delta} (j) \leq f_{\beta\delta} (j)$. Set $n = \sup \{j \in D_\delta : f_{\alpha\delta} (j) \not= f_{\beta\delta} (j)\}$. Since $\delta \notin S_{\alpha\beta}$, we have that $n < \delta$. On the other hand, $\delta \in S_\alpha$, so there is a cofinal subset $W$ of $\delta$ with $f_\alpha \vert W \subseteq A^{i}_\delta$. Now pick $j \in W$ with $j > \max (n, m)$. Then $(j, f_\alpha (j)) \in A^{i}_\delta$, and consequently $j \in D_\delta$. Hence $f_\alpha (j) = f_{\alpha\delta} (j) = f_{\beta\delta} (j) \geq f_\beta (j)$. This contradiction  completes the proof of the claim. 

\medskip

{\bf Claim 2.}  Let $\alpha < \beta < \gamma < \tau$. Then $\vert S_{\beta\gamma} \setminus S_{\alpha\gamma} \vert < \kappa$ and $\vert S_{\alpha\beta} \setminus S_{\alpha\gamma} \vert < \kappa$.

\medskip

{\bf Proof of Claim 2.}  Let $k_0$ (respectively, $k_1$) in $\kappa$ be such that $f_\alpha (j) < f_\beta (j)$ (respectively, $f_\beta (j) < f_\gamma (j)$) for all $j$ greater than $k_0$ (respectively, $k_1$). Then for any $j$ greater than $k_0$ (respectively, $k_1$), $f_{\alpha\xi} (j) \leq f_{\beta\xi} (j)$ (respectively,  $f_{\beta\xi} (j) \leq f_{\gamma\xi} (j)$) for all $\xi < \kappa$, and consequently $S_{\beta\gamma} \setminus S_{\alpha\gamma} \subseteq k_0$ (respectively,  $S_{\alpha\beta} \setminus S_{\alpha\gamma}
 \subseteq k_1$), which  completes the proof of the claim. 

\medskip

Since  there is no ascending $(J, I_\kappa)$-tower of length $\tau$, we may find, for each $\alpha < \tau$, $\alpha^\ast$ with $\alpha < \alpha^\ast < \tau$ such that $S_{\alpha\beta} \setminus S_{\alpha\alpha^\ast} \in J$ whenever
$\alpha^\ast < \beta < \tau$. 

\medskip

{\bf Claim 3.}  Let $\alpha < \beta < \tau$. Then $S_{\beta\beta^\ast} \setminus S_{\alpha\alpha^\ast} \in J$.

\medskip

{\bf Proof of Claim 3.}  Pick $\gamma$ with $\max (\alpha^\ast, \beta^\ast) < \gamma < \tau$. Then 

\centerline{$S_{\beta\beta^\ast} \setminus S_{\alpha\alpha^\ast} \subseteq (S_{\beta\beta^\ast} \setminus S_{\beta\gamma}) \cup (S_{\beta\gamma} \setminus S_{\alpha\gamma}) \cup (S_{\alpha\gamma} \setminus S_{\alpha\alpha^\ast})$,}

which  completes the proof of the claim. 

\medskip

Since  there is no descending $(J,J)$-tower of length $\tau$, we may find $\gamma < \tau$ such that $S_{\gamma\gamma^\ast} \setminus S_{\beta\beta^\ast} \in J$ whenever $\gamma < \beta < \tau$. Select $T \in J^+ \cap P (S_{\gamma\gamma^\ast})$ and $\theta < \nu$ such that $\vert A^{i}_\delta \vert = \theta$ for all $\delta \in T$. Inductively define $g : \theta^+ \rightarrow \tau \setminus (\gamma + 1)$ by : $g (\zeta)$ equals $\gamma + 1$ if $\zeta = 0$, and $(\sup \{g (\xi)^\ast : \xi < \zeta\}) + 1$ otherwise. Notice that if $\xi < \zeta < \theta^+$, then $\{S_{g (\xi)g (\xi)^\ast} \bigtriangleup S_{g (\xi)g (\zeta)}, S_{g (\xi)g (\xi)^\ast} \bigtriangleup S_{\gamma\gamma^\ast}\} \subseteq J$, so we may find $C_{\xi\zeta} \in J^\ast$ such that $S_{g (\xi)g (\zeta)} \cap C_{\xi\zeta} = S_{\gamma\gamma^\ast} \cap C_{\xi\zeta}$. Set $C = \bigcap \{C_{\xi\zeta} : \xi < \zeta < \theta^+\}$. Then $T \cap C \subseteq S_{g (\xi)g (\zeta)}$ whenever $\xi < \zeta < \theta^+$. Put 

\centerline{$s = (\sup \bigcup_{\xi < \zeta < \theta^+} \{j < \kappa : f_{g (\xi)} (j) \geq f_{g (\zeta)} (j) \}) + 1$,}

 and pick $\delta \in T \cap C$ with $\delta > s$. Notice that since $\delta \in T$, we have $\vert R_\delta \vert \leq \vert A^{i}_\delta \vert < \theta^+$. For each $j \in \kappa \setminus s$, the sequence $\langle f_{g (\xi)} (j) : \xi < \theta^+ \rangle$ is strictly increasing. It follows that for any $j \in D_\delta \setminus s$, the sequence  $\langle f_{g (\xi)\delta} (j) : \xi < \theta^+ \rangle$ is nondecreasing, and in fact eventually constant since $\{f_{g (\xi)\delta} (j) : \xi < \theta^+\} \subseteq R_\delta \cup \{\kappa\}$. Thus we may find $\xi_j <  \theta^+$ such that $f_{g (\xi)\delta} (j) = f_{g (\xi_j)\delta} (j)$ whenever $\xi_j < \xi < \theta^+$. Put $\eta = \sup \{\xi_j : j \in D_j \setminus s\}$, and let $\eta < \xi < \zeta < \theta^+$.  Then $f_{g (\xi)\delta} \vert (D_j \setminus s) = f_{g (\zeta)\delta} \vert (D_j \setminus s)$. However $\delta \in S_{g (\xi)g (\zeta)}$, and consequently $\ sup \{j \in D_\delta : f_{g (\xi)\delta} (j) < f_{g (\zeta)\delta} (j)\} = \delta$. Contradiction ! 
\hfill$\square$    

\bigskip

\subsection{Silly meeting}

\medskip

\begin{Def}  Given a regular infinite cardinal $\theta$, and a cardinal $\nu > \theta$, we let $M(\theta,\nu) =$  the least cardinality of any $X\subseteq P_\nu (\nu)$  with the property that  for any $e \in [\nu]^\theta$,  there is $x \in X$  with $\vert x \cap e \vert = \theta$.
\end{Def}

\begin{Obs} Let $\theta$ be a regular infinite cardinal, and $\nu > \theta$ be a cardinal with $\cf (\nu) \not= \theta$. Then $M(\theta, \nu) = \cf (\nu)$.
\end{Obs}

{\bf Proof.} Select $x_i \in P_\nu (\nu)$ for $i < \cf (\nu)$ such that
\begin{itemize}
\item $x_j \subseteq x_i$ for all $j < i$.
\item $\bigcup_{i < \cf (\nu)} x_i = \nu$.
\end{itemize}
Now given $e \in [\nu]^\theta$, define $f : e \rightarrow \cf (\nu)$ by $f (\alpha) =$ the least $j$ such that $\alpha \in x_j$. There must be $i < \cf (\nu)$ such that $\vert f^{-1} (i+ 1) \vert = \theta$. Then clearly $\vert x_i \cap e \vert = \theta$. 
\hfill$\square$


\begin{fact}  {\rm (\cite{Rin10})}  Let $\theta$ be a regular infinite cardinal less than $\kappa$. Suppose that $\kappa = \nu^+$, where $\cf (\nu )\not= \theta$. Then there is ${\cal B}_\delta \subseteq P_\nu (\delta))$  with $\vert {\cal B}_\delta \vert \leq \cf (\nu)$ for $\delta < \kappa$ such that for any $A \in [\kappa]^\kappa$,

\centerline{$\{ \delta < \kappa : \exists B \in {\cal B}_\delta (\sup (A \cap B) = \delta)\} \in (NS_\kappa \vert E^\kappa_\theta)^\ast$.}
\end{fact}

{\bf Proof.} Using Observation 3.11, for $\delta \in E^\kappa_\theta \setminus \nu$, pick ${\cal B}_\delta \subseteq P_\nu (\delta)$  with $\vert {\cal B}_\delta \vert \leq \cf (\nu)$ such that for any $e \in [\delta]^\theta$,  there is $B \in {\cal B}_\delta$  with $\vert B \cap e \vert = \theta$. Given $A \in [\kappa]^\kappa$, set $C = \{ \delta \in acc (\kappa) \setminus \nu : \sup (A \cap \delta) = \delta\}$. Now fix $\delta \in C \cap E^\kappa_\theta$. Select $e \subseteq A \cap \delta$ with $\sup e = \delta$ and o.t.$(e) = \theta$. We may find $B \in {\cal B}_\delta$ such that $\vert B \cap e \vert = \theta$. Then clearly, $\sup (B \cap A) = \delta$.  
\hfill$\square$

\begin{Rmk} It obviously follows that if $\kappa = \nu^+$, where $\cf (\nu )\not= \theta$, then $\clubsuit_\kappa^{{\rm cof}, - /Êcf (\nu)} [J]$ holds for any $\kappa$-complete ideal $J$ on $\kappa$ extending $NS_\kappa \vert E^\kappa_\theta$. For the case when $\theta = \cf (\nu) < \nu$ and $J = NS_\kappa \vert S$ for some $S \in NS^+_\kappa \cap P (E^\kappa_\theta)$, see Theorem 2.6 in \cite{Rin10} which gives a condition in terms of approachability for $\clubsuit_\kappa^{{\rm cof}, - /Êcf (\nu)} [J]$ to hold.
\end{Rmk}

\begin{Pro}  Let $\theta$ be a regular infinite cardinal less than $\kappa$. Suppose that $\kappa = \nu^+$, where $\cf (\nu )\not= \theta$, and $J$ is a $\kappa$-complete ideal on $\kappa$ extending $NS_\kappa \vert E^\kappa_\theta$. Then for any regular cardinal $\tau < Depth ({}^\kappa \kappa)$, there exists either an ascending $(J, I_\kappa)$-tower of length $\tau$, or a descending $(J, J)$-tower of length $\tau$. 
\end{Pro}

{\bf Proof.} By Theorem 3.9 and Remark 3.13.
\hfill$\square$

\bigskip

\subsection{Slow train}

\medskip

We will now give a proof of Fact 2.18. Let us recall the setting : $\theta < \kappa$ is a regular uncountable cardinal, $\kappa = \nu^+ = 2^\nu$, where $\cf (\nu) \not= \theta$, and $J$ is a $\kappa$-complete ideal on $\kappa$ extending $NS_\kappa \vert E^\kappa_\theta$. By Remark 3.13, we already know that $\clubsuit_\kappa^{{\rm cof}, - /Êcf (\nu)} [J]$ holds. We need to show that $\clubsuit_\kappa [J]$ holds. Just as Primavesi \cite{Prima}, we are not looking for a concise, beautiful proof. On the contrary, the more steps the better, as we would like to see in slow motion how $\clubsuit_\kappa^{{\rm cof}, - /Êcf (\nu)} [J]$ gradually evolves into $\clubsuit_\kappa [J]$. The main component of the proof is assertion (i) in the following proposition.

\medskip

\begin{Pro}  \begin{enumerate}[\rm (i)]
\item Suppose that $\clubsuit_\kappa^{{\rm cof} / \sigma, - /Ê\tau} [J]$ holds, where $\sigma$ is an infinite cardinal with $\kappa^\sigma = \kappa$, $\tau$ is a cardinal with $1 \leq \tau < \kappa$, and $J$ is a $\kappa$-complete ideal on $\kappa$. Then $\clubsuit_\kappa^{- /Ê\tau} [J]$ holds.
\item Suppose that $\clubsuit_\kappa^{{\rm cof} / \sigma, -} [J]$ holds, where $\sigma$ is an infinite cardinal with $\kappa^\sigma = \kappa$ and $J$ is a $\kappa$-complete ideal on $\kappa$. Then $\clubsuit_\kappa^- [J]$ holds.
\end{enumerate}
\end{Pro}

{\bf Proof.} (i) : We modify the proof (which we do not understand) of Theorem 3.5 in \cite{Prima} that asserts that assertion (i) is valid for any ideal of the form $NS_\kappa \vert S$.
Let $s^{i}_\gamma \in P_\sigma (\gamma)$ for $i < \tau$ and $\gamma \in acc (\kappa)$ witness that $\clubsuit_\kappa^{{\rm cof} / \sigma, - /Ê\tau} [J]$ holds. Let $\langle v_\delta : \delta < \kappa \rangle$ be a one-to-one enumeration of ${}^\sigma \kappa$. Define $Bad : \sigma \times [\kappa]^\kappa \rightarrow P(\kappa)$ by $Bad (r, A) = \{ \delta \in \kappa : v_\delta(r) \notin A\}$. Given $r < \sigma$ and $k : r \rightarrow [\kappa]^\kappa$, define $f_k : \tau \times acc (\kappa) \rightarrow P_\sigma (\kappa)$ by 
 
\centerline{$f_k (i, \gamma) = \{v_\delta (r) : \delta \in s^{i}_\gamma \setminus (\bigcup_{q < r} Bad(q, k(q)))\}$.}
  
\medskip

{\bf Claim.}   There is $k \in\bigcup_{r < \sigma}{}^r (P(\kappa))$ such that 

\centerline{$\{\gamma \in acc (\kappa) : \exists i < \tau (\sup f_k (i, \gamma) = \gamma \, {\rm and} \, f_k (i, \gamma) \subseteq T)\} \in J^+$} 

for all $T \in [\kappa]^\kappa$.

\medskip

{\bf Proof of the claim.}  Suppose otherwise.  For each $k \in\bigcup_{r < \sigma}{}^r ([\kappa]^\kappa)$, pick $T_k \in [\kappa]^\kappa$ and $C_k \in J^\ast \cap P (acc (\kappa))$ such that for any $\gamma \in C_k$ and any $i < \tau$, either $\sup f_k (i, \gamma) < \gamma$, or $f_{k} (i, \gamma) \setminus T_k  \not= \emptyset$. Define $H : \sigma \rightarrow [\kappa]^\kappa$ so that $H(r) = T_{H \vert r}$ for all $r < \sigma$. For $r < \sigma$, let $\langle \xi^r_\beta : \beta < \kappa \rangle$ be the increasing enumeration of $H (r)$. Define $F : \kappa \rightarrow \kappa$ so that $v_{F (\beta)}(r) = \xi^r_\beta$ for every $\beta < \kappa$ and every $r < \sigma$. Notice that $F$ is one-to-one. Inductively define $\beta_j < \kappa$ for $j < \kappa$ so that $\sup \{ \max (\beta_l, F (\beta_l)) : l < j\} < \beta_j$. Put $\Delta = \{F (\beta_j) : j < \kappa\}$. There must be $\gamma \in \bigcap_{r < \sigma} C_{H \vert r}$ and $i < \tau$ such that $\sup (s^{i}_\gamma \cap \Delta) = \gamma$.

Since $\vert s^{i}_\gamma \vert < \sigma$, we may find $r < \sigma$ such that 

\centerline{$(s^{i}_\gamma \setminus \bigcup_{q < r} Bad(q, H(q)) \cap Bad(r, H(r)) = \emptyset$.}

 Then $v_\delta(r) \in H(r)$ for any $\delta \in s^{i}_\gamma \setminus (\bigcup_{q < r} Bad(q, H(q))$, and therefore $f_{H|r}(i, \gamma) \subseteq H(r)$. Given $\alpha < \gamma$, pick $l < j < \kappa$ so that $\{F (\beta_l), F (\beta_j) \} \subseteq s^{i}_\gamma$ and $F (\beta_l) \geq \alpha$. Then $v_{F (\beta_j)} (r) = \xi^r_{\beta_j} \geq \beta_j > F (\beta_l) \geq \alpha$. Now clearly, $s^{i}_\gamma \cap \Delta \subseteq s^{i}_\gamma \setminus Bad(u, H(u))$ for all $u < \kappa$, so $s^{i}_\gamma \cap \Delta \subseteq s^{i}_\gamma \setminus (\bigcup_{q < r} Bad(q, H(q)))$. Hence, $v_{F (\beta_j)} (r) \in f_{H|r}(i, \gamma) \setminus \alpha$. Thus $\sup f_{H|r}(i, \gamma) = \gamma$. This contradiction completes the proof of the claim and that of (i).
 
 (ii) : The proof is a straightforward modification of that of (i).    
\hfill$\square$

\medskip

Primavesi \cite{Prima} established that if $\kappa^\tau = \kappa$ and $\clubsuit_\kappa^{- / \tau} [NS_\kappa \vert S]$ holds, where $S  \in NS_\kappa^+$, then so does $\clubsuit_\kappa [NS_\kappa \vert S]$. This can be generalized as follows.

\medskip

\begin{Pro}  
Let $\tau$ be a cardinal with $1 < \tau < \kappa$ and $\kappa^\tau = \kappa$, and $J$ be a $\kappa$-complete ideal on $\kappa$ extending $NS_\kappa$.  Suppose that $\clubsuit_\kappa^{- / \tau} [J]$ holds. Then so does $\clubsuit_\kappa [J]$.
\end{Pro}
 
{\bf Proof.} 
The proof is a modification of that of Theorem 6.2.3 in \cite{Prima}. Thus let  $s^{i}_\alpha \subseteq \alpha$ with $\sup s^{i}_\alpha = \alpha$ for $\alpha \in acc (\kappa)$ and $i < \tau$ be such that $\{ \alpha \in acc (\kappa) : \exists i < \tau (s^{i}_\alpha \subseteq A)\} \in J^+$ for every $A \in [\kappa]^\kappa$. Let $\langle e_\gamma : \gamma < \kappa \rangle$ be a $\kappa$-to-one enumeration of $[\tau \times \kappa]^\tau$. For $\alpha \in acc (\kappa)$ and $i < \tau$, put $t^{i}_\alpha = \bigcup_{\gamma \in s^{i}_\alpha} \{\delta < \kappa : (i, \delta) \in e_\gamma\}$. 

We claim that there is $i < \tau$ such that for any $A \in [\kappa]^\kappa$, the set of all $\alpha \in acc (\kappa)$ such that $\sup t^{i}_\alpha = \alpha$ and $t^{i}_\alpha \subseteq A$ lies in $J^+$. Suppose otherwise. Then we may find $C_i \in J^\ast \cap P( acc (\kappa))$ and $A_i \in [\kappa]^\kappa$ for $i < \tau$ such that for any $i < \tau$ and any $\alpha \in C_i$, either $\sup t^{i}_\alpha \not= \alpha$, or $t^{i}_\alpha \setminus A_i \not= \emptyset$. For $i < \tau$, let $\langle a^{i}_\beta : \beta < \kappa \rangle$ be the increasing enumeration of $A_i$. We inductively define $\beta_\xi, \gamma_\xi < \kappa$ for $\xi < \kappa$ as follows. We let $\beta_0 = \gamma_0 = 0$. Assuming that $\xi > 0$ and $\beta_\zeta$ and $\gamma_\zeta$ have been constructed for all $\zeta < \xi$, we let $\beta_\xi =$ the least $\beta > \sup \{\beta_\zeta : \zeta < \xi\}$ such that $\min \{a^{i}_\beta : i < \tau\} \geq \xi$, and $\gamma_\xi =$ the least $\gamma > \sup \{\gamma_\zeta : \zeta < \xi\}$ such that $e_\gamma = \{(i, a^{i}_{\beta_\xi}) : i < \tau\}$. Let $D$ be the set of all $\alpha \in acc (\kappa)$ such that
\begin{itemize}
\item $\gamma_\alpha = \alpha$ ;
\item $\alpha \in \bigcap_{i < \tau} C_i$ ;
\item $\{\delta < \kappa : (i, \delta) \in e_\gamma\} \subseteq \alpha$ for any $\gamma < \alpha$ and and any $i < \tau$.
\end{itemize}
There must be $\alpha \in D$ and $i < \tau$ such that  $s^{i}_\alpha \subseteq \{\gamma_\xi : \xi < \kappa\}$. It is readily checked that $t^{i}_\alpha \subseteq A_i \cap \alpha$. Given $i < \zeta < \alpha$, we may find $\xi < \kappa$ such that $\gamma_\xi \in s^{i}_\alpha \setminus \gamma_\zeta$. Since $i < \zeta \leq \xi \leq \beta_\xi$, we have $(i, a^{i}_{\beta_\xi}) \in e_{\beta_\xi}$, and therefore $a^{i}_{\beta_\xi} \in t^{i}_\alpha$. Furthermore, $\zeta \leq \xi \leq a^{i}_{\beta_\xi}$. Thus $\sup t^{i}_\alpha = \alpha$. Contradiction.
\hfill$\square$

\medskip

Assertion (i) in the following proposition is due to Rinot \cite{Rin10} in the case when $J$ is a restriction of $NS_\kappa$.

\medskip

\begin{Pro}  \begin{enumerate}[\rm (i)]
\item Suppose that $\kappa$ is a successor cardinal, and $J$ is a $\kappa$-complete ideal on $\kappa$ extending $NS_\kappa$. Then the following are equivalent :
\begin{enumerate}[\rm (a)]
\item $\diamondsuit_\kappa [J]$ holds.
\item $\clubsuit_\kappa^{{\rm cof}, -} [J]$ holds and $2^{< \kappa} = \kappa$.
\end{enumerate}
\item Suppose that $\kappa$ is weakly inaccessible, and $J$ is a normal ideal on $\kappa$. Then for any infinite cardinal $\sigma < \kappa$, the following are equivalent : 
\begin{enumerate}[\rm (i)]
\item $\diamondsuit_\kappa [J]$ holds.
\item $\clubsuit_\kappa^{{\rm cof} / \sigma, -} [J]$ holds and $2^{< \kappa} = \kappa$.
\end{enumerate}
\end{enumerate}
\end{Pro}

{\bf Proof.}  (i) : By Observation 2.22 and Propositions 3.15 and 3.16.

(ii) : By Fact 2.16, Observation 3.4 and Propositions 3.15 and 3.16.
\hfill$\square$

\medskip

{\bf Proof of Fact 2.18.} By Remark 3.13 and Proposition 3.17.
\hfill$\square$

\medskip

Mildenberger \cite{Heike} showed that for any $S \in NS^+_{\omega_1}$, if CH and $\clubsuit_\kappa^{\rm ev} [NS_{\omega_1} \vert S]$ both hold, then $\diamondsuit_{\omega_1} [NS_{\omega_1} \vert S]$ holds (see \cite{KS} for more results of this type). This generalizes.

\medskip

\begin{Obs}  Let $\theta$ be a regular infinite cardinal less than $\kappa$, and $J$ be a $\kappa$-complete ideal on $\kappa$ extending $NS_\kappa \vert E^\kappa_\theta$ such that $\clubsuit_\kappa^{\rm ev} [J]$ holds. Then the following hold : 
\begin{enumerate}[\rm (i)]
\item Suppose that $\kappa^\theta = \kappa$. Then  $\clubsuit_\kappa [J]$ holds.
\item Suppose that $2^{< \kappa} = \kappa$. Then $\diamondsuit_\kappa [J]$ holds. 
\end{enumerate}
\end{Obs}

{\bf Proof.} Observe that $\clubsuit_\kappa^{- / \theta} [J]$ holds, and appeal to Proposition 3.16 and Observation 2.22.
\hfill$\square$

\medskip

By considering other versions of the club principle, one can obtain variants of Proposition 3.16.

\medskip

\begin{Def}  Given a cardinal $\tau$ with $1 \leq \tau < \kappa$ and a $\kappa$-complete ideal $J$ on $\kappa$, $\clubsuit_\kappa^{\rm ev / -Ê\tau} [J]$ asserts the existence of  $s^{i}_\alpha \subseteq \alpha$ with $\sup s^{i}_\alpha = \alpha$ for $\alpha \in acc (\kappa)$ and $i < \tau$ such that $\bigcup_{i < \tau} \{ \alpha \in acc (\kappa) :  \exists \beta < \alpha (s^{i}_\alpha \setminus \beta \subseteq A\}$ lies in $J^+$ for all $A \in [\kappa]^\kappa$.
\end{Def}

\begin{Obs}  Let $\tau$ be a cardinal with $1 < \tau < \kappa$ and $\kappa^\tau = \kappa$, and $J$ be a $\kappa$-complete ideal on $\kappa$ extending $NS_\kappa$.  Then the following hold : 
\begin{enumerate}[\rm (i)]
\item Suppose that $\clubsuit_\kappa^{\rm ev / -Ê\tau} [J]$ holds. Then so does $\clubsuit_\kappa^{\rm ev} [J]$.
\item Suppose that there are $s^{i}_\alpha \subseteq \alpha$ for $\alpha \in acc (\kappa)$ and $i < \tau$ such that $\bigcup_{i < \tau} \{ \alpha \in acc (\kappa) : \sup (s^{i}_\alpha \cap A) = \alpha\} \in J^+$ for all $A \in [\kappa]^\kappa$. Then there are  $t_\alpha \subseteq \alpha$ with $\vert t_\alpha \vert \leq \max (\tau, \vert s^{i}_\alpha \vert)$ for $\alpha \in acc (\kappa)$ such that $ \{ \alpha \in acc (\kappa) : \sup (t_\alpha \cap A) = \alpha\} \in J^+$ for all $A \in [\kappa]^\kappa$.
 \end{enumerate}
\end{Obs}


\bigskip

\subsection{The case $\kappa = \nu^+$ with $\nu$ singular}

\medskip

We start by recalling the definition of covering numbers.

\medskip

\begin{Def}  Given four cardinals $\rho_1, \rho_2, \rho_3, \rho_4$ with $\rho_1 \geq \rho_2 \geq \rho_3 \geq \omega$ and $\rho_3 \geq \rho_4 \geq 2$, $\cov (\rho_1, \rho_2, \rho_3, \rho_4)$  denotes the least cardinality of any $Z \subseteq  P_{\rho_2}(\rho_1)$ such that for any $a \in  P_{\rho_3}(\rho_1)$, there is $Q \in  P_{\rho_4}(Z)$ with $a \subseteq \bigcup Q$.          
\end{Def}

\begin{Obs}  Let $\tau$, $\chi$ and $\sigma$ be three cardinals such that $1 \leq \tau \leq \chi < \sigma < \kappa$ and $\cov (\kappa, \chi^+, \tau^+, 2) = \kappa$. Suppose that $\clubsuit_\kappa^{{\rm cof} / \sigma, - /Ê\tau} [J]$ holds, where $J$ is a $\kappa$-complete ideal on $\kappa$. Then $\clubsuit_\kappa^{{\rm cof} / \sigma} [J]$ holds.
\end{Obs}

{\bf Proof.} The proof is similar to that of Proposition 3.16. Let $s^{i}_\alpha \in P_\sigma (\alpha)$ for $i < \tau$ and $\alpha \in acc (\kappa)$ witness that $\clubsuit_\kappa^{{\rm cof} / \sigma, - /Ê\tau} [J]$ holds. Pick $Z \subseteq  P_{\chi^+}(\kappa)$ such that $\vert Z \vert = \kappa$ and for any $a \in  P_{\tau^+}(\kappa)$, there is $z \in  Z$ with $a \subseteq z$.  Let $\langle z_\gamma : \gamma < \kappa \rangle$ be a $\kappa$-to-one enumeration of $Z$. For $\alpha \in acc (\kappa)$ and $i < \tau$, put $t^{i}_\alpha = \alpha \cap (\bigcup_{\gamma \in s^{i}_\alpha} z_\gamma)$. Note that $\vert t^{i}_\alpha \vert \leq \max \{ \vert s^{i}_\alpha \vert, \chi\} < \sigma$.

We claim that there is $i < \tau$ such that for any $A \in [\kappa]^\kappa$, the set of all $\alpha \in acc (\kappa)$ such that $\sup (t^{i}_\alpha \cap A) = \alpha$ lies in $J^+$. Suppose otherwise. Then we may find $C_i \in J^\ast \cap P( acc (\kappa))$ and $A_i \in [\kappa]^\kappa$ for $i < \tau$ such that for any $i < \tau$ and any $\alpha \in C_i$, $\sup (t^{i}_\alpha \cap A_i) < \alpha$. For $i < \tau$, let $\langle a^{i}_\beta : \beta < \kappa \rangle$ be the increasing enumeration of $A_i$. We inductively define $\beta_\xi, \gamma_\xi < \kappa$ for $\xi < \kappa$ as follows. We let $\beta_0 = \gamma_0 = 0$. Assuming that $\xi > 0$ and $\beta_\zeta$ and $\gamma_\zeta$ have been constructed for all $\zeta < \xi$, we let $\beta_\xi =$ the least $\beta > \sup \{\beta_\zeta : \zeta < \xi\}$ such that $\min \{a^{i}_\beta : i < \tau\} \geq \xi$, and $\gamma_\xi =$ the least $\gamma > \sup \{\gamma_\zeta : \zeta < \xi\}$ such that $\{a^{i}_{\beta_\xi} : i < \tau\} \subseteq z_\gamma$. Let $D$ be the set of all $\alpha \in acc (\kappa) \cap \bigcap_{i < \tau} C_i$ such that $\bigcup_{\gamma < \alpha} z_\gamma \subseteq \alpha = \gamma_\alpha$. There must be $\alpha \in D$ and $i < \tau$ such that $\sup (s^{i}_\alpha \cap \{\gamma_\xi : \xi < \kappa\}) = \alpha$. Given $i < \zeta < \alpha$, we may find $\xi < \kappa$ such that $\gamma_\xi \in s^{i}_\alpha \setminus \gamma_\zeta$. Then clearly $\xi < \alpha$ and $z_{\gamma_\xi} \subseteq \alpha$. Since $a^{i}_{\beta_\xi} \in z_{\gamma_\xi}$, it follows that $a^{i}_{\beta_\xi} \in t^{i}_\alpha$. Furthermore, $\zeta \leq \xi \leq a^{i}_{\beta_\xi}$. Thus $\sup (t^{i}_\alpha \cap A_i) = \alpha$. Contradiction.
\hfill$\square$

\begin{fact} Suppose that $\kappa = \nu^+$, where $\nu$ is singular. Then 

\centerline{$\cov (\nu, \nu, (\cf (\nu))^+, 2) = \cov (\kappa, \chi^+, (\cf (\nu))^+, 2)$}

 for some cardinal $\chi < \nu$.
\end{fact}
 
{\bf Proof.}  Set $\theta = \cov (\nu, \nu, (\cf (\nu))^+, 2)$. By \cite[Observation 5.3 (10) p. 86]{SheCA}, there must be a cardinal $\chi < \nu$ such that $\theta = \cov (\nu, \chi, (\cf (\nu))^+, 2)$. Then clearly, $\theta = \cov (\nu, \chi^+, (\cf (\nu))^+, 2)$. On the other hand, by \cite[Observation 5.3 (2) p. 86]{SheCA}, $\cov (\kappa, \chi^+, (\cf (\nu))^+, 2) = \max \{\theta, \kappa\}$. It remains to observe that by \cite[Theorem 5.4 p. 88]{SheCA}, $\theta \geq \pp (\nu)$.
\hfill$\square$

\begin{Pro}  Let $\theta$ be a regular infinite cardinal less than $\kappa$. Suppose that $\kappa = \nu^+$, where $\cf (\nu )\in \nu \setminus \{\theta\}$ and $\cov (\nu, \nu, (\cf (\nu))^+, 2) = \kappa$, and $J$ is a $\kappa$-complete ideal on $\kappa$ extending $NS_\kappa \vert E^\kappa_\theta$. Then $\clubsuit_\kappa^{{\rm cof}} [J]$ holds. 
\end{Pro} 
 
{\bf Proof.} $\clubsuit_\kappa^{{\rm cof}, - /Ê\tau} [J]$ holds by Remark 2.10, and hence so does $\clubsuit_\kappa^{{\rm cof}} [J]$ by Observation 3.22 and Fact 3.23.
 \hfill$\square$

\medskip

Note that by results of Shelah (\cite[Theorem 5.4 p. 87]{SheCA}, \cite{She93}), SSH is equivalent to the statement that $\cov (\nu, \nu, (\cf (\nu))^+, 2) = \nu^+$ for any singular cardinal $\nu$.

\bigskip

\subsection{Slow train II}

\medskip

This time we would like to retrace the path leading from Fact 3.12 to Fact 2.17.

\medskip

\begin{Def}  Given a $\kappa$-complete ideal $J$ on $\kappa$ and a cardinal $\sigma < \kappa$, $\clubsuit^{\rm cof / \sigma, \ast}_\kappa [J]$ asserts the existence of $B^{i}_\delta \in P_\sigma (\delta))$  for $i < \delta < \kappa$ such that for any $A \in [\kappa]^\kappa$,

\centerline{$\{ \delta < \kappa : \exists i < \delta (\sup (A \cap B^{i}_\delta) = \delta)\} \in J^\ast$.}
\end{Def}

\medskip

We will follow Fuchs and Rinot who established \cite{Fuchs} that if $\kappa = \nu^+ = 2^\nu$ and $\clubsuit^{\rm cof / \nu, \ast}_\kappa [NS_\kappa \vert S]$ holds for $S \in NS^+_\kappa \vert E^\kappa_\theta$, where $\theta$ is a regular infinite cardinal less than $\nu$ such that $d (\theta, \mu) \leq \nu$ for any cardinal $\mu$ with $\theta \leq \mu < \nu$, then $\diamondsuit^\ast_\kappa [NS_\kappa \vert S]$ holds. However we will make an extra stop at $\clubsuit_\kappa^\ast$.

\medskip

\begin{Pro}  Let $\theta < \sigma$ be two infinite cardinals such that $\cf (\theta) = \theta$, $\sigma < \kappa$, and $d (\theta, \mu) < \kappa$ for any cardinal $\mu$ with $\theta \leq \mu < \sigma$. Further let $J$ be a $\kappa$-complete ideal on $\kappa$ extending $NS_\kappa \vert E^\kappa_\theta$ such that $\clubsuit_\kappa^{\rm cof / \sigma, \ast} [J]$ holds. Then $\clubsuit_\kappa^\ast [J]$ holds. 
\end{Pro}

{\bf Proof.} Let $B^{i}_\delta \in P_\sigma (\delta))$  for $i < \delta < \kappa$ witness that $\clubsuit_\kappa^{\rm cof / \sigma, \ast} [J]$ holds. For $i < \delta < \kappa$, we define $Z^{i}_\delta$ as follows. If $\vert B^{i}_\delta \vert < \theta$, put $Z^{i}_\delta = \emptyset$. Otherwise let $Z^{i}_\delta$ be a cofinal subset of $([B^{i}_\delta]^\theta, \supseteq)$ of size $d (\theta, \vert B^{i}_\delta \vert)$. For $\delta < \kappa$, set ${\cal A}_\delta = \bigcup \{ Z^{i}_\delta : i < \delta$ \, {\rm and} \,$ \vert B^{i}_\delta \vert \geq \theta \}$. Note that if $\delta \geq \sup \{d (\theta, \mu) : \theta \leq \mu < \sigma\}$, then $\vert {\cal A}_\delta \vert \leq \vert \delta \vert$. 
Given $A \in [\kappa]^\kappa$, we may find $C \in J^\ast$ such that for any $\delta \in C$, there is $i < \delta$ with $\sup (A \cap B^{i}_\delta) = \delta$. Now fix $\delta \in C \cap E^\kappa_\theta$.  Pick $i < \delta$ and $e \subseteq A \cap B^{i}_\delta$ so that $\sup e = \delta$ and o.t.$(e) = \theta$. There must be $z \in Z^{i}_\delta$ such that $z \subseteq e$. Then $z \in {\cal A}_\delta$. Moreover, $z \subseteq A$ and $\sup z = \delta$. 
\hfill$\square$

\begin{Cor} Suppose that $\kappa = \nu^+$, and $\theta$ is a regular infinite cardinal less than $\nu$ such that $d (\theta, \nu) = \nu$. Then $\clubsuit^\ast_\kappa [NS_\kappa \vert E^\kappa_\theta]$ holds.  
\end{Cor}

{\bf Proof.} Use Fact 3.12.
\hfill$\square$

\medskip

{\bf Proof of Fact 2.17.} By Observation 2.23 and Corollary 3.27.
\hfill$\square$

\bigskip

\subsection{More clubbing}

\medskip

Corollary 3.27 can be easily generalized to weakly inaccessible cardinals, but it is not yet clear what these weak generalizations (already considered in \cite{DZ0} and \cite {She10}) are good for. The following strengthens a result of D\v zamonja \cite{DZ0}.

\medskip

\begin{Obs} There is ${\cal B}_\delta \subseteq \{B \subseteq \delta : \sup B = \delta\}$ with $\vert {\cal B}_\delta \vert \leq d (\cf (\delta), \vert \delta \vert)$  for $\delta \in acc (\kappa)$ such that for any $A \in [\kappa]^\kappa$, there exists $C \in {\cal C}_\kappa$ with the property that $C \subseteq \{ \delta \in acc (\kappa) : \exists B \in {\cal B}_\delta (B \subseteq A)\}$. 
\end{Obs}

{\bf Proof.} For $\delta \in acc (\kappa)$, pick ${\cal A}_\delta \subseteq [\delta]^{\cf (\delta)}$ with $\vert {\cal A}_\delta \vert = d (\cf (\delta), \vert \delta \vert)$ such that for any $e \in [\delta]^{\cf (\delta)}$, there is $B \in {\cal A}_\delta$ with $B \subseteq e$. Given $A \in [\kappa]^\kappa$, let $C = \{ \delta \in acc (\kappa) : \sup (A \cap \delta) = \delta \}$. Now fix $\delta \in C$. Select $e \subseteq A \cap \delta$ such that o.t.$(e) = \cf (\delta)$ and $\sup e = \delta$. There must be $B \in {\cal A}_\delta$ such that $B \subseteq e$. Then clearly $\sup B = \delta$, and moreover $B \subseteq A$.

\hfill$\square$

\bigskip

\subsection{Order-type versus cardinality}

\medskip

To conclude this section let us mention the following result of D\v zamonja and Shelah that yields a variant of $\clubsuit_\kappa^{{\rm cof}}$ where the condition on $\vert s_\alpha \vert$ is replaced with one on o.t.$(s_\alpha)$.   

\medskip

\begin{fact}  {\rm (\cite{DZ3})}  Suppose that $\kappa = \nu^+$, where $\nu$ is regular, and $\theta$ and $\rho$ are two regular infinite cardinals with $\theta < \rho < \nu$. Suppose further that either $\theta > \omega$, or $\nu \geq 2^{\aleph_0}$, and let $S \in NS_\kappa^+ \cap P (E^\kappa_\theta)$ with the property that $S$ reflects at stationarily many $\delta \in E^\kappa_\rho$. Then there is $s_\alpha \subseteq \alpha$ for $\alpha \in S$ with {\rm o.t.}$(s_\alpha) < \nu^{+ \omega} \cdot \rho$ such that $\{\alpha \in S : \sup (A \cap s_\alpha) = \alpha\} \in NS^+_\kappa$ for all $A \in [\kappa]^\kappa$. 
\end{fact}

\bigskip

\section{GUESSING GENERALIZED CLUBS}
\bigskip

In this section we revisit another result of Rinot where nonsaturation is derived from guessing generalized clubs. Let us start with some definitions.

\medskip

\begin{Def}  Let $\sigma$ be an infinite cardinal, and $\delta$ be a limit ordinal greater than or equal to $\sigma$. A subset $C$ of $P_\sigma (\delta)$ is a {\it generalized club} if $\{x \in P_\sigma (\delta) : F`` P_\omega (x) \subseteq x\} \subseteq C$ for some $F : P_\omega (\delta) \rightarrow \delta$.
\end{Def}

\begin{Obs} \begin{enumerate}[\rm (i)]
\item If $\sigma = \omega$, then there is an empty generalized club subset of $P_\sigma (\delta)$.
\item If $\sigma > \omega$, then any generalized club subset of $P_\sigma (\delta)$ is cofinal in $(P_\sigma (\delta), \subseteq)$.
\end{enumerate}
 \end{Obs}

\begin{Def}  Given an infinite cardinal $\sigma < \kappa$ and a $\kappa$-complete ideal $J$ on $\kappa$, $\curlywedge^- (\sigma, \kappa, J)$) asserts the existence of a cofinal subset $C^{i}_\delta$ of $(P_\sigma (\delta), \subseteq)$ for $i \in \delta \in acc (\kappa) \setminus \sigma$ such that $\{\delta : \exists i < \delta \, (C^{i}_\delta \subseteq D)\} \in J^+$ for every generalized club subset $D$ of $P_\sigma (\kappa)$. 

$\curlywedge (\sigma, \kappa, J)$ asserts the existence of a generalized club subset $C_\delta$ of $P_\sigma (\delta)$ for $\delta \in acc (\kappa) \setminus \sigma$ such that $\{\delta : C_\delta \subseteq D\} \in J^+$ for every generalized club subset $D$ of $P_\sigma (\kappa)$. 
\end{Def}

\medskip

Note that $\curlywedge (\sigma, \kappa, J) \Rightarrow \curlywedge^- (\sigma, \kappa, J) \Rightarrow \sigma > \omega$. $\curlywedge (\sigma, \kappa, NS_\kappa \vert S)$ (respectively, $\curlywedge^- (\sigma, \kappa, NS_\kappa \vert S)$) is identical with the principle $\curlywedge (\sigma, S)$ (respectively, $\curlywedge^- (\sigma, S)$) introduced in \cite{RinGen}. For the associated starred version see \cite{KLY}. 

\begin{Obs} The following are equivalent :
\begin{enumerate}[\rm (i)]
\item $\curlywedge^- (\sigma, \kappa, J)$ holds.
\item There is a subset $X^{i}_\delta$ of $P_\sigma (\delta)$ with the property that $\delta \subseteq \bigcup X^{i}_\delta$ for $i \in \delta \in acc (\kappa) \setminus \sigma$ such that $\{\delta : \exists i < \delta \, (X^{i}_\delta \subseteq D)\} \in J^+$ for every generalized club subset $D$ of $P_\sigma (\kappa)$. 
\end{enumerate}
\ \end{Obs}

{\bf Proof.} It suffices to prove that (ii) implies (i) since the other direction is trivial. Thus let $X^{i}_\delta$ for $i \in \delta \in acc (\kappa) \setminus \sigma$ be as in (ii). For $\delta \in acc (\kappa) \setminus \sigma$, select a cofinal subset $Z_\delta$ of $P_\sigma (\delta)$, and $h : \delta \times \delta \rightarrow X^{i}_\delta$ such that $\alpha \in h (i, \alpha)$ for each $(i, \alpha) \in \delta \times \delta$. Set $C^{i}_\delta = \{\bigcup_{\alpha \in z} h (i, \alpha) : z \in Z_\delta\}$. Now given a generalized club subset $D$ of $P_\sigma (\kappa)$, pick $F : P_\omega (\kappa) \rightarrow \kappa$ such that $D' \subseteq D$, where $D' = \{x \in P_\rho (\delta) : F`` P_\omega (x) \subseteq x\}$. Then clearly, $C^{i}_\delta \subseteq D' \subseteq D$ for any $\delta \in acc (\kappa) \setminus \sigma$ with $X^{i}_\delta \subseteq D'$. 
\hfill$\square$

\begin{Cor} If $\curlywedge^- (\sigma, \kappa, J)$ holds, then so does $\curlywedge^- (\sigma', \kappa, J)$ for every cardinal $\sigma'$ with $\sigma \leq \sigma' < \kappa$.
\ \end{Cor}

\medskip

The following is essentially due to Shelah.

\medskip

\begin{Obs} Suppose that $J$ extends $NS_\kappa$, and $\diamondsuit_\kappa [J]$ holds. Then for any ucountable cardinal $\sigma < \kappa$, $\curlywedge (\sigma, \kappa, J)$ holds.
 \end{Obs}

{\bf Proof.} Let $\langle s_\delta : \delta < \kappa \rangle$ witness that $\diamondsuit_\kappa [J]$ holds. Select a bijection $j : P_\omega (\kappa) \times \kappa \rightarrow \kappa$. For $\delta \in acc (\kappa)$, define $f_\delta : P_\omega (\delta) \rightarrow \delta$ by : $f_\delta (e)$ equals  $0$ if $\{\xi \in \delta : j (e, \xi) \in s_\gamma\} = \emptyset$, and $\min \{\xi \in \delta : j (e, \xi) \in s_\gamma\}$ otherwise. Now let $\sigma$ be a fixed uncountable cardinal less than $\kappa$. Set $C_\delta  = \{x \in P_\sigma (\delta) : f_\delta`` P_\omega (x) \subseteq x\}$. Given a generalized club subset $D$ of $P_\sigma (\kappa)$ and $H \in J^\ast$, pick $F : P_\omega (\kappa) \rightarrow \kappa$ with  $\{x \in P_\omega (\kappa) : F`` P_\omega (x) \subseteq x\} \subseteq D$. We may find $\delta \in H$ such that $F`` P_\omega (\delta) \subseteq \delta = j`` (P_\omega (\delta) \times \delta)$ and $s_\delta = \{j (e, F (e)) : e \in P_\omega (\kappa)\} \cap \delta$. Then it is readily checked that $C_\delta \subseteq D$.
\hfill$\square$

\medskip

Rinot \cite{RinGen} established that if $\kappa = \nu^+$, where $\nu$ is regular, and $S \in NS^+_\kappa \cap P (E^\kappa_\nu)$ is such that  $\curlywedge^- (\sigma, \kappa, NS_\kappa \vert S)$ holds for some $\sigma < \kappa$, then $NS_\kappa \vert S$ is not $\kappa^+$-saturated. This can be extended as follows.

\medskip

\begin{Th}   \begin{enumerate}[\rm (i)]
\item Suppose that 
\begin{itemize}
\item $\kappa = \nu^+$ ;
\item $\sigma$ is an uncountable cardinal less than $\nu$ ;
\item $J$ is a $\kappa$-complete ideal on $\kappa$ such that $E^\kappa_{> \sigma} \in J^\ast$ and $\curlywedge^- (\sigma, \kappa, J)$ holds ; 
\item $\tau$ is a regular cardinal less than $Depth ({}^\kappa \kappa)$.
\end{itemize}
Then there exists either an ascending $(J, I_\kappa)$-tower of length $\tau$, or a descending $(J, J)$-tower of length $\tau$.
\item Suppose that 
\begin{itemize}
\item $\kappa$ is weakly inaccessible ;
\item $\sigma$ is an uncountable cardinal less than $\kappa$ ;
\item $J$ is a normal ideal on $\kappa$ such that $E^\kappa_{> \sigma} \in J^\ast$ and $\curlywedge^- (\sigma, \kappa, J)$ holds ; 
\item $\tau$ is a regular cardinal less than $Depth ({}^\kappa \kappa)$.
\end{itemize}
Then there exists either an ascending $(J, I_\kappa)$-tower of length $\tau$, or a descending $(J, J)$-tower of length $\tau$.

\end{enumerate}
\end{Th}

{\bf Proof.} We prove (ii) and leave the similar proof of (i) to the reader. The proof is a modification of that of Theorem 3.4 in \cite{RinGen}. Select $f_\alpha \in {}^\kappa \kappa$ for $\alpha < \tau$ such that $f_\alpha <^\ast f_\beta$ whenever $\alpha < \beta < \tau$. Let $C^{i}_\delta$ for $i \in \delta \in acc (\kappa) \setminus \sigma$ witness that $\curlywedge^- (\sigma, \kappa, J)$ holds. 
For $\delta \in E^\kappa_{> \sigma}$, select $h_\delta : \delta \times \delta \rightarrow P_\sigma (\delta)$ such that $j \in h_\delta (i, j) \in C^{i}_\delta$ for all $(i, j) \in \delta \times \delta$. 
For $\delta \in E^\kappa_{> \sigma}$ and $\alpha < \tau$, define $f^\delta_\alpha : \delta \times \delta \rightarrow \delta \cup \{\kappa\}$ by $f_\alpha^\delta (i, j) = \min ((h_\delta (i, j) \cup \{\kappa\}) \setminus f_\alpha (j))$. Notice that if $f_\alpha (j) \leq f_\beta (j)$, then $f_\alpha^\delta (i, j) \leq f_\beta^\delta (i, j)$. Finally, for $\alpha < \beta < \tau$ and $i < \kappa$, let

\centerline{$S^{i}_{\alpha\beta} = \{\delta \in E^\kappa_{> \sigma} \setminus (i + 1) :  \sup \{j \in \delta : f^\delta_\alpha (i, j) = f_\beta^\delta (i, j)\} < \delta\}.$}

\medskip

{\bf Claim 1.}  Let $\alpha < \beta < \gamma < \tau$ and $i < \kappa$. Then the following hold :
\begin{enumerate}[\rm (i)]
\item $\vert S^{i}_{\alpha\beta} \setminus S^{i}_{\alpha\gamma} \vert < \kappa$.
\item $\vert S^{i}_{\beta\gamma} \setminus S^{i}_{\alpha\gamma} \vert < \kappa$.
\end{enumerate}
\medskip

{\bf Proof of Claim 1.}  Pick $\eta < \kappa$ so that $f_\alpha (j) < f_\beta (j) < f_\gamma (j)$ whenever $\eta \leq j < \kappa$. Let us first show that $S^{i}_{\alpha\beta} \setminus (\eta + 1) \subseteq S^{i}_{\alpha\gamma}$. Given $\delta \in S^{i}_{\alpha\beta} \setminus (\eta + 1)$, put $A = \{j < \delta : f^\delta_\alpha (i, j) = f^\delta_\beta (i, j)\} \cup \eta$. Then clearly, $f^\delta_\alpha (i, j) < f^\delta_\beta (i, j) \leq f^\delta_\gamma (i, j)$ for all $j \in \delta \setminus A$. Thus $\{j < \delta : f^\delta_\alpha (i, j) = f^\delta_\gamma (i, j)\} \subseteq A$, and consequently $\delta \in S^{i}_{\alpha\gamma}$. Next we show that $S^{i}_{\beta\gamma} \setminus (\eta + 1) \subseteq S^{i}_{\alpha\gamma}$.  Given $\delta \in S^{i}_{\beta\gamma} \setminus (\eta + 1)$, put $B = \{j < \delta : f^\delta_\beta (i, j) = f^\delta_\gamma (i, j)\} \cup \eta$. Then clearly, $f^\delta_\alpha (i, j) \leq f^\delta_\beta (i, j) < f^\delta_\gamma (i, j)$ for all $j \in \delta \setminus B$. Hence $\{j < \delta : f^\delta_\alpha (i, j) = f^\delta_\gamma (i, j)\} \subseteq B$, and therefore $\delta \in S^{i}_{\alpha\gamma}$, which  completes the proof of the claim. 

\medskip

{\bf Claim 2.}  There is $v : \tau \rightarrow \kappa$ such that $S^{v (\alpha)}_{\alpha\beta} \in J^+$ whenever $\alpha < \beta < \tau$.

\medskip

{\bf Proof of Claim 2.} Fix $\alpha < \kappa$.  Set $k = (\sup \{j < \kappa : f_\alpha (j) \geq f_{\alpha + 1} (j)\}) + 1$. Define $F : P_\omega (\kappa) \rightarrow \kappa$ by :
\begin{itemize}
\item $F (e) = 0$ if $e = \emptyset$ ;
\item $F (e) = \xi + 1$ if $e = \{\xi\}$ :
\item $F (e) = f_\alpha( \max e)$ if $\vert e \vert = 2$ ;
\item $F (e) = f_{\alpha + 1} (\max e)$ if $\vert e \vert > 2$.
\end{itemize}
Put $D = \{x \in P_\sigma (\kappa) : F`` P_\omega (x) \subseteq x\}$. Note that $\{f_\alpha (\zeta), f_{\alpha + 1} (\zeta)\} \subseteq x$ whenever $x \in D$ and $\zeta \in x \setminus \omega$. By normality of $J$, there must be $i < \kappa$ such that 

\centerline{$T = \{\delta \in E^\kappa_{> \sigma} \setminus \max (i + 1, k + 1) : C^{i}_\delta \subseteq D\}$}

lies  in $J^+$. If $\delta \in T$, then for any $j \in \delta$ with $j \geq \max (k, \omega)$, we have that $j \in h_\delta (i, j) \in D$, so $f^\delta_\alpha (i, j) = f_\alpha (j) < f_{\alpha + 1} (j) = f^\delta_{\alpha + 1} (i, j)$. Thus $T \subseteq S^{i}_{\alpha\alpha+1}$. Hence $S^{i}_{\alpha(\alpha+1)}$ lies in $J^+$, and by Claim 1 so does $S^{i}_{\alpha\beta}$ for every $\beta > \alpha$, which  completes the proof of the claim. 

\medskip

There must be $i < \kappa$ such that $\vert v^{- 1} (\{i\}) \vert = \tau$. By thinning out our sequence of functions, we may assume that $v (\alpha) = i$ for all $\alpha < \tau$. 

Since  there is no ascending $(J, I_\kappa)$-tower of length $\tau$, we may find, for each $\alpha < \tau$, $\alpha^\ast$ with $\alpha < \alpha^\ast < \tau$ such that $S^{i}_{\alpha\beta} \setminus S^{i}_{\alpha\alpha^\ast} \in J$ whenever $\alpha^\ast < \beta < \tau$. 

\medskip

{\bf Claim 3.}  Let $\alpha < \beta < \tau$. Then $S^{i}_{\beta\beta^\ast} \setminus S^{i}_{\alpha\alpha^\ast} \in J$.

\medskip

{\bf Proof of Claim 3.}  Pick $\gamma$ with $\max (\alpha^\ast, \beta^\ast) < \gamma < \tau$. Then 

\centerline{$S^{i}_{\beta\beta^\ast} \setminus S{i}_{\alpha\alpha^\ast} \subseteq (S^{i}_{\beta\beta^\ast} \setminus S^{i}_{\beta\gamma}) \cup (S^{i}_{\beta\gamma} \setminus S^{i}_{\alpha\gamma}) \cup (S^{i}_{\alpha\gamma} \setminus S^{i}_{\alpha\alpha^\ast})$,}

which  completes the proof of the claim. 

\medskip

Since  there is no descending $(J,J)$-tower of length $\tau$, we may find $\gamma < \tau$ such that $S_{\gamma\gamma^\ast} \setminus S_{\beta\beta^\ast} \in J$ whenever $\gamma < \beta < \tau$. Select $T \in J^+ \cap P (S^{i}_{\gamma\gamma^\ast})$ and $\theta < \sigma$ such that $ \sup \{j < \delta : \vert h_\delta (i, j) \vert = \theta \} = \delta$ for all $\delta \in T$. Inductively define $g : \theta^+ \rightarrow \tau \setminus (\gamma + 1)$ by : $g (\zeta)$ equals $\gamma + 1$ if $\zeta = 0$, and $(\sup \{g (\xi)^\ast : \xi < \zeta\}) + 1$ otherwise. 

Notice that if $\xi < \zeta < \theta^+$, then $\{S^{i}_{g (\xi)g (\xi)^\ast} \bigtriangleup S^{i}_{g (\xi)g (\zeta)}, S^{i}_{g (\xi)g (\xi)^\ast} \bigtriangleup S^{i}_{\gamma\gamma^\ast}\} \subseteq J$, so we may find $C_{\xi\zeta} \in J^\ast$ such that $S^{i}_{g (\xi)g (\zeta)} \cap C_{\xi\zeta} = S^{i}_{\gamma\gamma^\ast} \cap C_{\xi\zeta}$. Set $C = \bigcap \{C_{\xi\zeta} : \xi < \zeta < \theta^+\}$. Then $T \cap C \subseteq S^{i}_{g (\xi)g (\zeta)}$ whenever $\xi < \zeta < \theta^+$. Put  

\centerline{$s = (\sup \bigcup_{\xi < \zeta < \theta^+} \{j < \kappa : f_{g (\xi)} (j) \geq f_{g (\zeta)} (j) \}) + 1$,}

 and pick $\delta \in T \cap C$ with $\delta > s$. For $\xi < \theta^+$, set $W_\xi = \{j < \delta : f_{g (\xi} (j) > \sup h_\delta (i, j)\}$.
 
 \medskip

{\bf Claim 4.}  Let $\xi < \theta^+$. Then $\sup W_\xi < \delta$.

\medskip

{\bf Proof of Claim 4.}  Suppose otherwise. Pick $\zeta$ with $\xi < \zeta < \theta^+$. Then $\kappa = f_{g (\xi}^\delta (i, j) \leq  f_{g (\zeta}^\delta (i, j) \leq \kappa$, contradicting the fact that $\delta \in S^{i}_{g (\xi)g (\zeta)}$, which  completes the proof of the claim. 

\medskip

Set $t = (\sup \bigcup_{\xi < \theta^+} W_\xi) + 1$, $u = \max (s, t)$ and $Q = \{j \in \delta \setminus u : \vert h_\delta (i, j) \vert = \theta\}$. For $j \in Q$, $\langle f^\delta_{g (\xi)} (i, j) : \xi < \theta^+ \rangle$ is a weakly increasing sequence of elements of $h_\delta (i, j)$, since the sequence $\langle f_{g (\xi)} (j) : \xi < \theta^+ \rangle$ is increasing, so there must be $\chi_j < \theta^+$ such that $f^\delta_{g (\xi)} (i, j) = f^\delta_{g (\chi_j)} (i, j)$ whenever $\chi_j < \xi < \theta^+$. We may find $\chi < \theta^+$ and $M \subseteq Q$ with $\sup M = \delta$ such that $\chi_j = \chi$ for all $j \in M$. But now $f^\delta_{g (\xi)} \vert (\{i\} \times M) = f^\delta_{g (\xi)} \vert (\{i\} \times M)$ whenever $\chi < \xi < \zeta < \theta^+$. Contradiction ! 
\hfill$\square$

\medskip

Let us now compare our principle for guessing generalized clubs with the principles considered in the previous section (a weak version of club) and the next section (club-guessing). The following extends Theorem 2.1 of \cite{RinGen}. 

\medskip

\begin{Obs} \begin{enumerate}[\rm (i)]
\item Suppose that $\kappa$ is weakly inaccessible and $\curlywedge^- (\sigma, \kappa, J)$ holds, where $J$ is a normal ideal on $\kappa$ such that $E^\kappa_{\geq \sigma} \in J^\ast$. Then there is a closed unbounded subset $X_\delta$ of $\delta$ for $\delta \in E^\kappa_{\geq \sigma}$ such that $\{\delta : \sup (X_\delta \setminus Y) < \delta\} \in J^+$ for every closed unbounded subset $Y$ of $\kappa$.
\item Suppose that $\kappa$ is a successor cardinal and $\curlywedge^- (\sigma, \kappa, J)$ holds, where $J$ is a $\kappa$-complete ideal on $\kappa$ such that $E^\kappa_{\geq \sigma} \in J^\ast$. Then there is a closed unbounded subset $X_\delta$ of $\delta$ for $\delta \in E^\kappa_{\geq \sigma}$ such that $\{\delta : X_\delta \subseteq Y\} \in J^+$ for every closed unbounded subset $Y$ of $\kappa$.
\end{enumerate}
\ \end{Obs}

{\bf Proof.} We prove (i) and leave the similar proof of (ii) to the reader. Let $C^{i}_\delta$ for $i \in \delta \in acc (\kappa) \setminus \sigma$ witness that $\curlywedge^- (\sigma, \kappa, J)$ holds. For $i \in \delta \in E^\kappa_{\geq \sigma}$, put $S^{i}_\delta = \{ \sup x : x \in C^{i}_\delta \setminus \{\emptyset\}\}$. Note that $S^{i}_\delta \subseteq \delta$.  

 \medskip

{\bf Claim 1.}  Let $Y$ be a closed unbounded subset of $\kappa$. Then $\{\delta \in E^\kappa_{\geq \sigma} : \exists i < \delta \, (S^{i}_\delta \subseteq Y)\} \in J^+$.

\medskip

{\bf Proof of Claim 1.} Define $F : P_\omega (\kappa) \rightarrow \kappa$ by $F (e) = \min (C \setminus \bigcup e)$, and let $D = \{x \in P_\sigma (\delta) : F`` P_\omega (x) \subseteq x\}$. Now suppose that $i$ and $\delta$ are such that $i \in \delta \in E^\kappa_{\geq \sigma}$ and $C^{i}_\delta \subseteq D$. Then clearly, $\sup x \in Y$ for every nonempty $x$ in $C^{i}_\delta$. Hence $S^{i}_\delta \subseteq Y$, which  completes the proof of the claim. 

\medskip
 For  $i \in \delta \in E^\kappa_{\geq \sigma}$, put $T^{i}_\delta = \{ \xi \in \delta \setminus \{0\} : \sup (\xi \cap T^{i}_\delta) = \xi\}$. Then clearly, $T^{i}_\delta$ is a closed unbounded subset of $\delta$. Furthermore, $\{\delta \in E^\kappa_{\geq \sigma} : \exists i < \delta \, (T^{i}_\delta \subseteq Y)\} \in J^+$ for any closed unbounded subset $Y$ of $\kappa$.
 \medskip

{\bf Claim 2.}  There is $i < \kappa$ such that  $\{\delta : T^{i}_\delta \setminus Y\subseteq i + 1\} \in J^+$ for every closed unbounded subset $Y$ of $\kappa$.

\medskip

{\bf Proof of Claim 2.}  Suppose otherwise. For $i < \kappa$, pick $Y_i \in NS_\kappa^\ast$ and $W_i \in J^\ast \cap P (E^\kappa_{\geq \sigma})$ such that $(T^{i}_\delta \setminus Y_i) \setminus (i + 1) \not= \emptyset$ whenever $i \in \delta \in W_i$. Set $Y = \bigtriangleup_{i < \kappa} Y_i$ and $W = \bigtriangleup_{i < \kappa} W_i$. We may find $\delta \in W$ and $i < \delta$ such that $T^{i}_\delta \subseteq Y$. Then $\delta \in W_i$, and moreover $Y^{i}_\delta \setminus (i + 1) \subseteq Y_i \setminus (i + 1) \subseteq Y$. This contradiction  completes the proof of the claim and that of the observation. 

\hfill$\square$

\begin{Obs} Suppose that $\curlywedge^- (\sigma, \kappa, J)$ holds, where $J$ is a $\kappa$-complete ideal on $\kappa$ such that $E^\kappa_{\geq \sigma} \in J^\ast$. Then $NS^\ast_\kappa \cap J = \emptyset$.
\ \end{Obs}

{\bf Proof.} By the proof of Observation 4.8.
\hfill$\square$

\begin{Obs} Suppose that $\curlywedge^- (\sigma, \kappa, J)$ holds and $E^\kappa_\theta \in J^\ast$, where $\theta$ is an infinite cardinal with $\theta^+ < \kappa$. Then $\clubsuit_\kappa^{{\rm cof} / \rho, -} [J]$ holds, where $\rho = \max \{\sigma, \theta^+\}$.
 \end{Obs}

{\bf Proof.} Let $C^{i}_\delta$ for $i \in \delta \in acc (\kappa) \setminus \sigma$ witness that $\curlywedge^- (\sigma, \kappa, J)$ holds. For $\delta \in E^\kappa_\theta \setminus \sigma$, pick an increasing sequence $\langle \delta_j : j < \theta \rangle$ with supremum $\delta$, and $h : \delta \times \theta \rightarrow \delta$ such that $\delta_j \in h (i, j) \in C^{i}_\delta$ for every $(i, j) \in \delta \times \theta$.  Set $s^{i}_\delta = \bigcup_{j < \theta} h (i, j)$ for every $i < \delta$. Now given $A \in [\kappa]^\kappa$, define $F : P_\omega (\kappa) \rightarrow \kappa$ by $F (e) = \min (A \setminus \cup e)$, and put $D = \{x \in P_\sigma (\delta) : F`` P_\omega (x) \subseteq x\}$. Then clearly, $\sup (A \cap s^{i}_\delta) = \delta$ whenever $\delta \in E^\kappa_\theta \setminus \sigma$ and $i < \delta$ are such that $C^{i}_\delta \subseteq D$. 
\hfill$\square$

\medskip

The following is yet another indication of the strength of $\curlywedge^-$.

\medskip

\begin{Obs} Given a regular uncountable cardinal $\sigma < \kappa$, the following hold :
\begin{enumerate}[\rm (i)]
\item Suppose that $\curlywedge (\sigma, \kappa, NS_\kappa \vert S)$ holds, where $S \in NS_\kappa^+ \cap P (E^\kappa_{< \sigma})$. Then there is a stationary subset $X$ of $P_\sigma (\kappa)$ such that
\begin{itemize}
\item $\{\sup x : x \in X\} = S$ ;
\item the sup-function is one-to-one on $X$.
\end{itemize}
\item Suppose that $\curlywedge^- (\sigma, \kappa, NS_\kappa \vert S)$ holds, where $S \in NS_\kappa^+ \cap P (E^\kappa_{< \sigma})$. Then there is a stationary subset $X$ of $P_\sigma (\kappa)$ such that
\begin{itemize}
\item $\{\sup x : x \in X\} = S$ ;
\item $\vert \{x \in X : \sup x = \alpha \} \vert \leq \vert \alpha \vert$ for all $\alpha \in S$.
\end{itemize}
\end{enumerate}
\end{Obs}

{\bf Proof.} We prove (i) and leave the similar proof of (ii) to the reader. Let $\langle C_\delta : \delta \in acc (\kappa) \setminus \sigma \rangle$ witness that $\curlywedge (\sigma, \kappa, NS_\kappa \vert S)$ holds. We define $x_\delta$ for $\delta \in S$ as follows. Given $\delta \in S$, select $e_\delta \subseteq \delta$ so that o.t.$(e_\delta) = \cf (\delta)$ and $\sup e_\delta = \delta$. Inductively define $x^n_\delta \in C_\delta$ for $n < \omega$ so that
\begin{itemize}
\item $x^0_\delta = e_\delta$ ;
\item $x^n_\delta \cup \sup (x^n_\delta \cap \kappa) \subseteq x^{n + 1}_\delta$.
\end{itemize}
Finally, $X = \{\bigcup_{n < \omega} x^n_\delta : \delta \in S\}$ is as desired.
\hfill$\square$

\begin{Cor} Suppose that $\curlywedge^- (\sigma, \kappa, NS_\kappa \vert S)$ holds, where $\sigma$ is a regular uncountable cardinal less than $\kappa$, and $S \in NS_\kappa^+ \cap P (E^\kappa_{< \sigma})$. Then  $\cov (\kappa, \sigma, \sigma, 2) = \kappa$.
 \end{Cor}

\medskip

Notice that in the other direction, the following is known.

\medskip
  
\begin{fact}  {\rm (\cite{GS}, \cite{MSU})}  Let $S \in NS^+_\kappa \cap P (E^\kappa_{< \sigma})$. Suppose that there is a stationary subset $X$ of $P_\sigma (\kappa)$ such that
\begin{itemize}
\item $\{\sup x : x \in X\} \subseteq S$ ;
\item the sup-function is one-to-one on $X$ (respectively, $\vert \{x \in X : \sup x = \alpha \} \vert \leq \vert \alpha \vert$ for all $\alpha \in S$).
\end{itemize}
Then $\clubsuit_\kappa^{{\rm cof} / \sigma} [NS_\kappa \vert S]$ (respectively, $\clubsuit_\kappa^{{\rm cof} / \sigma, -} [NS_\kappa \vert S]$) holds. 
\end{fact}

\bigskip

\section{GITIK-SHELAH ON NONSATURATION}
\bigskip

We are looking for another result on ideal nonsaturation where towers might be involved. A natural candidate is the following result of Gitik and Shelah.

\medskip

\begin{fact} {\rm (\cite{GS})} Suppose that $\max (\omega_2, \theta^+) < \kappa$, where $\theta$ is a regular infinite cardinal. Then $NS_\kappa \vert E^\kappa_\theta$ is not $\kappa^+$-saturated.
 \end{fact}
 
 \begin{Rmk} Notice that in general the conclusion will not remain valid if $NS_\kappa \vert E^\kappa_\theta$ is replaced with $NS_\kappa \vert S$ for some stationary subset $S$ of $E^\kappa_\theta$ (see \cite{GS}).
 \end{Rmk}

\medskip

The result was already revisited by Krueger \cite{Kru}, and we will follow his reading.

\medskip

\subsection{Club-guessing }

\medskip

We start with an easy generalization of Shelah's club guessing principle (see \cite{GS}).

{\bf Throughout Subsections 5.1 - 5.4, $\theta$ and $\rho$ will denote two regular infinite cardinals with $\theta < \rho < \kappa$.}

\begin{Pro}  Let $K$ be a $\kappa$-complete ideal on $\kappa$ extending $NS_\kappa \vert E^\kappa_\theta$. Then there is $c_\alpha\subseteq E^\kappa_{\geq \rho} \cap \alpha$ with $\sup c_\alpha = \alpha$ for $\alpha \in E^\kappa_\theta \cap acc (E^\kappa_{\geq \rho})$ such that $\{\alpha \in E^\kappa_\theta \cap acc (E^\kappa_{\geq \rho}) :  c_\alpha \subseteq C\} \in K^+$ for every $C \in {\cal C}_\kappa$.  
\end{Pro}

{\bf Proof.} The proof follows that of Hirata given in \cite{Shi}. For each $\beta \in acc (\kappa)$, select a cofinal subset $d_\beta$ of $\beta$ of order-type $\cf (\beta)$. Given $\alpha \in E^\kappa_\theta \cap acc (E^\kappa_{\geq \rho})$ and $D \in {\cal C}_\kappa \cap P (acc (\kappa))$, we inductively define $x^D_{\alpha, n} \subseteq D$ with $\vert x^D_{\alpha, n} \vert < \rho$ for $n < \omega$ as follows :
\begin{itemize}
\item $x^D_{\alpha, 0} = \{ \sup (D \cap \gamma) : \gamma \in d_\alpha\}$ ;
\item $x^D_{\alpha, n + 1} = \{ \sup (D \cap \gamma) : \exists \beta \in x^D_{\alpha, n} \cap E^\kappa_{< \rho} (\gamma \in d_\beta)\}$.
\end{itemize}
Finally, we let $x^D_\alpha = (\bigcup_{n < \omega} x^D_{\alpha, n}) \setminus \{0\}$.

\medskip

{\bf Claim.}  There is $D \in {\cal C}_\kappa \cap P (acc (\kappa))$ such that $\{\alpha \in E^\kappa_\theta \cap acc (E^\kappa_{\geq \rho}) :  x^D_\alpha \subseteq C\} \in K^+$ for every $C \in {\cal C}_\kappa$. 

\medskip

{\bf Proof of the claim.}  Suppose otherwise. Inductively define $C_\xi \in {\cal C}_\kappa \cap P (acc (\kappa))$ and $H_\xi \in K^\ast$ for $\xi < \rho$ such that
\begin{itemize}
\item $C_{\xi + 1} \subseteq C_\xi$ ;
\item $H_{\xi + 1} \cap \{\alpha \in E^\kappa_\theta \cap acc (E^\kappa_{\geq \rho}) :  x^{C_\xi}_\alpha \subseteq C_{\xi + 1}\} = \emptyset$.
\end{itemize}
Select $\alpha$ in $E^\kappa_\theta \cap acc (E^\kappa_{\geq \rho})  \cap \bigcap_{\xi < \rho} C_\xi \cap \bigcap_{\xi < \rho} H_\xi$. Since the map $\xi \rightarrow \sup (C_\xi \cap \gamma)$ is nonincreasing for every $\gamma < \kappa$, and $\vert x^{C_\xi}_{\alpha, n} \vert < \rho$ for all $\xi < \rho$ and all $n < \omega$, we may inductively define $\xi_n$ for $n < \omega$ so that
\begin{itemize}
\item $\xi_n \leq \xi_{n + 1}$ ;
\item $x^{C_\xi}_{\alpha, n} = x^{C_{\xi_n}}_{\alpha, n}$ for $\xi_n \leq \xi < \rho$.
\end{itemize}
Set $\xi = \sup \{\xi_n : n < \omega \}$. Then $x^{C_\xi}_\alpha = x^{C_{\xi + 1}}_\alpha \subseteq C_{\xi + 1}$. This contradiction completes the proof of the claim.

\medskip

By the claim, for any $C \in {\cal C}_\kappa$,

\centerline{$\{\alpha \in E^\kappa_\theta \cap acc (E^\kappa_{\geq \rho}) \cap acc (D) :  x^D_\alpha \setminus acc (x^D_\alpha) \subseteq C \cap acc (D) \} \in K^+$.}

So it remains to observe the following. Let $\alpha \in E^\kappa_\theta \cap acc (E^\kappa_{\geq \rho}) \cap acc (D)$ be such that $x^D_\alpha \subseteq acc (D)$. Then $x^D_{\alpha, 0}$ is cofinal in $\alpha$, and hence so are $x^D_\alpha$ and $x^D_\alpha \setminus acc (x^D_\alpha)$. Furthermore for any $n < \omega$ and any $\beta \in x^D_{\alpha, n} \cap E^\kappa_{< \rho}$,  $\beta$ lies in $acc (D)$ and therefore in $acc (x^D_{\alpha, n + 1})$. It follows that $x^D_\alpha \setminus acc (x^D_\alpha) \subseteq E^\kappa_{\geq \rho}$.
\hfill$\square$  

\medskip

\subsection{Strong guessing }

\medskip

Let $J$ be a $\kappa$-complete ideal on $\kappa$ extending $NS_\kappa \vert E^\kappa_\theta$. By Proposition 5.3, for any $A \in J^+ \cap P (E^\kappa_\theta \cap acc (E^\kappa_{\geq \rho}))$, we may find $c^{A, J}_\alpha \subseteq E^\kappa_{\geq \rho} \cap \alpha$ for $\alpha \in A$ such that 
\begin{itemize}
\item for any $\alpha \in A$, $\sup c^{A, J}_\alpha = \alpha$ and o.t.$(c^{A, J}_\alpha) = \theta$ ;
\item $\{\alpha \in A :  c^{A,J}_\alpha \subseteq C\} \in J^+$ for every $C \in {\cal C}_\kappa$.
\end{itemize}

\medskip

\begin{Pro}  Let $J$ be a $\kappa$-complete ideal on $\kappa$ extending $NS_\kappa \vert E^\kappa_\theta$. Suppose that there is no descending $(J, I_\kappa)$-tower of length ${\frak b}_\kappa$. Then for any $A \in J^+ \cap P (E^\kappa_\theta \cap acc (E^\kappa_{\geq \rho}))$, there is $B \in J^+ \cap P (A)$ such that

\centerline{$\{\alpha \in A :  \exists \beta < \alpha (c^{A, J}_\alpha \setminus \beta \subseteq C)\} \in (J \vert B)^\ast$}

 for every $C \in {\cal C}_\kappa$.
\end{Pro}

{\bf Proof.} We follow the proof of Claim 2.1 in \cite{GS}. Fix $A \in J^+ \cap P (E^\kappa_\theta \cap acc (E^\kappa_{\geq \rho}))$. Let $\Phi : {\cal C}_\kappa \rightarrow J^+ \cap P (A)$ be defined by $\Phi (C) = \{\alpha \in A :  \exists \beta < \alpha (c^{A, J}_\alpha \setminus \beta \subseteq C)\}$. Notice that
\begin{itemize}
\item $\Phi (D) \subseteq \Phi (C)$ for all $D \in {\cal C}_\kappa \cap P (C)$ ;
\item $ \vert \Phi (D) \setminus \Phi (C) \vert < \kappa$ for all $D \in {\cal C}_\kappa$ such that $ \vert D \setminus C \vert < \kappa$.
\end{itemize}
Now suppose toward a contradiction that for any $B \in J^+ \cap P (A)$, there is $C \in {\cal C}_\kappa$ such that $B \setminus \Phi (C) \in J^+$. We inductively define $C_i \in {\cal C}_\kappa$ for $i < {\frak b}_\kappa$ so that for $j < i < {\frak b}_\kappa$, $\Phi (C_j) \setminus \Phi (C_i) \in J^+$ and $ \vert \Phi (C_i) \setminus \Phi (C_j) \vert < \kappa$. We put $C_0 = \kappa$. Now suppose that $i > 0$, and $C_j$ has been constructed for every $j < i$. By Observation 2.36 (i), we may find $D \in {\cal C}_\kappa$ such that $\vert D \setminus C_j \vert < \kappa$ for all $j < i$. By assumption, there must be $H \in {\cal C}_\kappa$ such that $\Phi (D) \setminus \Phi (H) \in J^+$. We let $C_i = D \cap H$. Finally, $\langle \Phi (C_i) : i < {\frak b}_\kappa \rangle$ is a descending $(J, I_\kappa)$-tower of length ${\frak b}_\kappa$, which yields the desired contradiction.
\hfill$\square$  
\bigskip

\subsection{Club-guessing ideals}

\medskip

\begin{Def}  Given a $\kappa$-complete ideal $J$ on $\kappa$ extending $NS_\kappa \vert E^\kappa_\theta$, we let ${\cal X}_J$ denote the collection of all $B \subseteq \kappa$ such that either $B \in J$, or $B \in J^+$ and there is $d_\alpha \subseteq E^\kappa_{\geq \rho} \cap \alpha$ for $\alpha \in B \cap E^\kappa_\theta \cap acc (E^\kappa_{\geq \rho})$ such that 
\begin{itemize}
\item for any $\alpha \in B \cap E^\kappa_\theta \cap acc (E^\kappa_{\geq \rho})$, $\sup d_\alpha = \alpha$ and o.t.$(d_\alpha) = \theta$ ;
\item $\{\alpha \in B \cap E^\kappa_\theta \cap acc (E^\kappa_{\geq \rho}) :  \exists \beta < \alpha (d_\alpha \setminus \beta \subseteq C)\} \in (J\vert B)^\ast$ for every $C \in {\cal C}_\kappa$.
\end{itemize}
\end{Def}

\begin{Obs}  \begin{enumerate}[\rm (i)]
\item Let $B \in {\cal X}_J$. Then $P (B) \subseteq {\cal X}_J$.
\item Let $A, B \in {\cal X}_J$. Then $A \cup B \in {\cal X}_J$.
\end{enumerate}
\end{Obs}

\begin{Rmk} Suppose that there is no descending $(J, I_\kappa)$-tower of length ${\frak b}_\kappa$. Then by Proposition 5.4, $ {\cal X}_J \cap P (H) \cap J^+ \not= \emptyset$ for any $H \in J^+$. 
\end{Rmk}

\begin{Pro} Let $J$ be a normal ideal on $\kappa$ with $E^\kappa_\theta \in J^\ast$. Suppose that there is no descending $(J, I_\kappa)$-tower of length ${\frak b}_\kappa$, and no ascending $(J, I_\kappa)$-tower of length $\kappa^+$. Then $\kappa \in {\cal X}_J$. 
\end{Pro}

{\bf Proof.} The proof is a modification of that of Lemma 2 in \cite{GS}. Let ${\cal T}$ denote the collection of those ascending $(J, I_\kappa)$-towers that have all their members in ${\cal X}_J$. Given $T, W \in {\cal T}$, put $T < W$ just in case $W$ is a proper extension of $T$. By Zorn's Lemma, $({\cal T}, <)$ has a maximal element, say $T = \langle B_\eta : \eta < \delta \rangle$. 

\medskip

{\bf Case 1 :} $\delta$ is a successor ordinal, say $\delta = \xi + 1$. 

\medskip

{\bf Claim 1.}  $\kappa \setminus B_\xi \in J$.

\medskip

{\bf Proof of Claim 1.}  Suppose otherwise. Then by Remark 5.7, we may find $R$ in ${\cal X}_J$ with $R \in J^+ \cap P (\kappa \setminus B_\xi)$. Set $B_\delta = B_\xi \cup R$ and $W = \langle B_\gamma : \gamma \leq \delta \rangle$. Then by Observation 5.6, $W \in {\cal T}$, and moreover $T < W$. This contradiction completes the proof of the claim.

\medskip

It follows from the claim and Observation 5.6 that $\kappa \in {\cal X}_J$.

\medskip

{\bf Case 2 :} $\delta$ is a limit ordinal. 

Put $\sigma = \cf (\delta)$, and let $\langle \delta_i : i < \sigma \rangle$ be an increasing sequence of ordinals with supremum $\delta$. For $i < \sigma$, set $K_i = B_{\delta_{i + 1}}$. Then it is simple to see that $\langle K_i : i < \sigma \rangle$ is also a maximal element of ${\cal T}$. Note that $\sigma \leq \kappa$. Set $H_0 = K_0 \setminus 1$, and for each $i$ with $0 < i < \sigma$, $H_i = K_i \setminus ((\bigcup_{j < i} K_j) \cup (i + 1))$. Notice that $H_m \cap H_n = \emptyset$ whenever $m < n < \sigma$. 

\medskip

{\bf Claim 2.}  Let $i < \sigma$. Then $H_i \in J^+$.

\medskip

{\bf Proof of Claim 2.}  This is obvious for $i = 0$. For $i > 0$, it suffices to observe that 

\centerline{$(B_{\delta_{i + 1}} \setminus B_{\delta_i}) \setminus \bigcup_{j < i} (B_{\delta_{j + 1}} \setminus B_{\delta_i}) \subseteq H_i$,}

which completes the proof of the claim.

\medskip

For each $i < \sigma$, $H_i \in {\cal X}_J$ by Observation 5.6, so we may find $d^{i}_\alpha \subseteq E^\kappa_{\geq \rho} \cap \alpha$ for $\alpha \in H_i \cap E^\kappa_\theta \cap acc (E^\kappa_{\geq \rho})$ such that 
\begin{itemize}
\item for any $\alpha \in H_i \cap E^\kappa_\theta \cap acc (E^\kappa_{\geq \rho})$, $\sup d^{i}_\alpha = \alpha$ and o.t.$(d^{i}_\alpha) = \theta$ ;
\item $\{\alpha \in H_i \cap E^\kappa_\theta \cap acc (E^\kappa_{\geq \rho}) :  \exists \beta < \alpha (d^{i}_\alpha \setminus \beta \subseteq C)\} \in (J \vert H_i)^\ast$ for every $C \in {\cal C}_\kappa$.
\end{itemize}
Put $A = \bigcup_{i < \sigma} H_i$.

\medskip

{\bf Claim 3.}  $A \in {\cal X}_J$.

\medskip

{\bf Proof of Claim 3.} For $\alpha \in A$, set $c_\alpha = d^{i}_\alpha$, where $\alpha \in H_i$. Now let $C \in {\cal C}_\kappa$. For each $i < \sigma$, we may find $D_i \in J^\ast$ such that for any $\alpha \in D_i \cap H_i \cap E^\kappa_\theta \cap acc (E^\kappa_{\geq \rho})$, there is $\beta < \alpha$ such that $c_\alpha \setminus \beta \subseteq C$. Put $D = \{ \alpha < \kappa : \forall i \in \sigma \cap \alpha (\alpha \in D_i)\}$. Notice that $D \in J^\ast$ by normality of $J$. Clearly, $H_i \cap D \subseteq D_i$ for all $i < \sigma$. It follows that for any $\alpha \in D \cap A \cap E^\kappa_\theta \cap acc (E^\kappa_{\geq \rho})$, there is $\beta < \alpha$ such that $c_\alpha \setminus \beta \subseteq C$, which completes the proof of the claim.

\medskip

{\bf Claim 4.}  $\kappa \setminus A \in J$.

\medskip

{\bf Proof of Claim 4.}  Suppose not. We proceed as in the proof of Claim 1. By Remark 5.7, there must be $R \in {\cal X}_J$ such that $R \in J^+ \cap P (\kappa \setminus A)$. Set $K_\delta = A \cup R$ and $W = \langle K_i : i \leq \sigma \rangle$. Now for each $i \leq \sigma$, $K_i \setminus A \subseteq i + 1$, and therefore $\vert K_i \setminus A \vert < \kappa$. By Observation 5.6, it follows that $W \in {\cal T}$. However, $T < W$. This contradiction completes the proof of the claim.

\medskip

By Claim 4 and Observation 5.6, $\kappa \in {\cal X}_J$.
\hfill$\square$  

\bigskip

\subsection{Full reflection}

\medskip

\begin{Def}  Given a stationary subset $S$ of $\kappa$, and a stationary subset $T$ of $E^\kappa_{\geq \omega_1}$, $S$ {\it reflects fully in} $T$ if there is $G \in {\cal C}_\kappa$ such that $S$ reflects at every $\gamma \in G \cap T$.
 \end{Def}
 
\begin{Obs} Let $\sigma < \kappa$ be a regular uncountable cardinal, and $S$ be a stationary subset of $\kappa$ that reflects fully in $E^\kappa_\sigma$. Then for any regular cardinal $\tau$ with $\sigma \leq \tau < \kappa$, $S$ reflects fully in $E^\kappa_\tau$.
\end{Obs}

\begin{fact} {\rm (\cite{Kru})} Suppose that $\kappa$ is either weakly inaccessible, or the successor of a singular cardinal. Then for any regular uncountable cardinal $\chi < \kappa$, and any stationary subset $S$ of $\kappa$, $S$ reflects fully in $E^\kappa_{\geq \chi}$ if and only if there are $C \in {\cal C}_\kappa$ and $\eta < \kappa$ such that $C \setminus S$ has no closed subset of order-type $\eta$.
\end{fact}

\begin{Th} Let $S$ be a stationary subset of $E^\kappa_\theta$ that reflects fully in $E^\kappa_{\geq \rho}$. Suppose that either $\theta > \omega$, or $\rho \geq \omega_2$. Then, setting $J = NS_\kappa \vert S$, there is either a descending $(J, I_\kappa)$-tower of length ${\frak b}_\kappa$, or an ascending $(J, I_\kappa)$-tower of length $\kappa^+$. 
\end{Th}

{\bf Proof.} We closely follow the proof of Lemma 3 in \cite{GS}. Suppose that the conclusion fails. Then by Proposition 5.4, there is $c_\alpha \subseteq E^\kappa_{\geq \rho} \cap \alpha$ for $\alpha \in S \cap acc (E^\kappa_{\geq \rho})$ such that 
\begin{itemize}
\item for any $\alpha \in S \cap acc (E^\kappa_{\geq \rho})$, $\sup c_\alpha = \alpha$ and o.t.$(c_\alpha) = \theta$ ;
\item $\{\alpha \in S \cap acc (E^\kappa_{\geq \rho}) :  \exists \beta < \alpha (c_\alpha \setminus \beta \subseteq C)\} \in J^\ast$ for any $C \in {\cal C}_\kappa$.
\end{itemize}
Pick $G \in {\cal C}_\kappa$ such that $S$ reflects at every $\gamma \in G \cap E^\kappa_{\geq \rho}$.
\medskip

{\bf Case 1 :} $\theta > \omega$. 

\medskip

Inductively define $C_n \in {\cal C}_\kappa$ for $n < \omega$ as follows. Set $C_0 = G \cap acc (E^\kappa_{\geq \rho})$. Now suppose that $C_n$ has been constructed. There must be $H \in {\cal C}_\kappa$ with the property that for any $\alpha \in H \cap S  \cap acc (E^\kappa_{\geq \rho})$, there is $\beta < \alpha$ such that $c_\alpha \setminus \beta \subseteq acc (C_n)$. Put $C_{n + 1} = acc (C_n) \cap H$. Finally, set $C = \bigcap_{n < \omega} C_n$ and $\alpha = \min (C \cap S)$. Since $\alpha \in \bigcap_{n < \omega} C_{n + 1}$, we may find $\eta < \alpha$ such that $c_\alpha \setminus \eta \subseteq \bigcap_{n < \omega} acc (C_n)$. Pick $\gamma \in c_\alpha \cap \bigcap_{n < \omega} acc (C_n)$. Then $\gamma \in E^\kappa_{\geq \rho} $. Moreover, $C_n \cap \gamma$ is a closed unbounded subset of $\gamma$ for every $n < \omega$. Hence $C \cap \gamma$ is a closed unbounded subset of $\gamma$, and therefore $(C \cap S) \cap \gamma \not= \emptyset$, which contradicts the minimality of $\alpha$.

\medskip

{\bf Case 2 :} $\theta = \omega$. 

In the same spirit as in Case 1, we define $C_\xi \in {\cal C}_\kappa$ for $\xi < \omega_1$ so that
\begin{itemize}
\item $C_0 =  G \cap acc (E^\kappa_{\geq \rho})$ ;
\item $C_\zeta \subseteq C_\xi$ for all $\zeta < \xi$ ;
\item for any $\alpha \in C_{\xi + 1} \cap S$, there is $\beta < \alpha$ such that $c_\alpha \setminus \beta \subseteq acc (C_\xi)$.
\end{itemize}
Set $C = \bigcap_{\xi < \omega_1} C_\xi$ and $\alpha = \min (C \cap S)$. Since $\alpha \in \bigcap_{\xi < \omega_1} C_{\xi + 1}$ and o.t.$(c_\alpha) = \omega$, there must be $\eta < \alpha$ such that $c_\alpha \setminus \eta \subseteq \bigcap_{\xi < \omega_1} acc (C_\xi)$. Pick $\gamma \in c_\alpha \cap \bigcap_{\xi < \omega_1} acc (C_\xi)$. Then $\cf (\gamma) \geq \rho > \omega_1$, and moreover $C \cap \gamma$ is a closed unbounded subset of $\gamma$. Hence $(C \cap S) \cap \gamma \not= \emptyset$. Contradiction.
\hfill$\square$   


\medskip

Full reflection is a sufficient condition, but in general not a necessary one, as the following shows.

\medskip

\begin{fact} {\rm (\cite{Kru})} Suppose that one of the following conditions holds :
\begin{itemize}
\item $\kappa$ is  weakly inaccessible.
\item  $\kappa = \nu^+$, where $\nu$ is singular. 
\item $\kappa = \nu^+$, where $\nu$ is regular and $\square_\nu$ holds.
 \end{itemize}
Then any stationary subset $T$ of $\kappa$ has a stationary subset $S$ with the property that for every regular uncountable cardinal $\sigma < \kappa$, $S$ does not reflect fully in $E^\kappa_\sigma$.
\end{fact}

\bigskip

\subsection{Good points}

Let $A$ be an infinite set of regular cardinals such that $\vert A \vert < \min A$ and $\sup A < \kappa$, and $I$ be an ideal on $A$ such that $\{A \cap a : a \in A \} \subseteq I$. 

\begin{Def}  We let $\prod A = \prod_{a \in A} a$. For $f, g \in \prod A$, we let $f <_I g$ if $\{a \in A : f (a) \geq g (a) \} \in I$. 
\end{Def}

\begin{Def} Let $\vec{f} = \langle f_\alpha : \alpha < \kappa \rangle$ be an increasing, cofinal sequence in $(\prod A, <_I)$. An infinite limit ordinal $\delta < \kappa$  is a {\it good point} for $\vec f$ if there is a cofinal subset $X \subseteq \delta$, and $Z_\xi \in I$ for $\xi \in X$ such that $f_\beta (a) < f_\xi (a)$ whenever $\beta < \xi$ are in $X$ and $a \in A \setminus (Z_\beta \cup Z_\xi)$.

We let $G (\vec{f})$ denote the set of good points for $\vec{f}$.  
\end{Def}

\begin{fact}  \begin{enumerate}[\rm (i)]
\item  {\rm(Folklore)}  If $\delta \in G (\vec{f})$, then $\cf(\delta) < \sup A$.
\item  {\rm(\cite{CFM}, \cite{Norm})} Let $\delta < \pi$  be an infinite limit ordinal such that $I$ is $\cf (\delta)$-complete. Then $\delta \in G (\vec{f})$.
\end{enumerate}
\end{fact}

\medskip

The following is due to Shelah.

\medskip

\begin{fact}  For $i = 0, 1$, let $\vec{f_i} = \langle f^{i}_\alpha : \alpha < \kappa \rangle$ be an increasing, cofinal sequence in $(\prod A, <_I)$. Then $G (\vec{f_0}) \bigtriangleup G (\vec{f_1}) \in NS_\kappa$.
\end{fact}

{\bf Proof.} Let $D$ be the set of all $\delta \in acc (\kappa)$ with the property that for any $\xi < \delta$, there are $\beta, \gamma < \delta$ such that $f^0_\xi <_I f^1_\beta$ and $f^1_\xi <_I f^0_\gamma$. Let us show that $D \cap G (\vec{f_0}) = D \cap G (\vec{f_1})$. Thus fix $i < 1$ and $\delta \in D \cap G (\vec{f_i})$. Let  $X \subseteq \delta$ and $Z_\xi \in I$ for $\xi \in X$ witness that $\delta$ is a good point for $\vec{f_i}$. Define two increasing sequences $\langle \beta_j : j < \cf (\delta) \rangle$ and $\langle \gamma_j : j < \cf (\delta) \rangle$ so that
\begin{itemize}
\item $\beta_j \in X$ and $\gamma_j < \delta$ ;
\item $f^{i}_{\beta_j} <_I f^{1 - i}_{\gamma_j} < _I f^{i}_{\beta_{j + 1}}$.
\end{itemize}
For $j < \cf (\delta)$, set

\centerline{$W_j = \{a : f^{i}_{\beta_j} (a) \geq f^{1 - i}_{\gamma_j} (a)\} \cup Z_{\beta_j} \cup \{a : f^{1 - i}_{\gamma_j} (a) \geq f^{i}_{\beta_{j + 1}} (a)\} \cup Z_{\beta_{j + 1}}$.} 

Then 

\centerline{$f^{1 - i}_{\gamma_j} (a) < f^{i}_{\beta_{j + 1}} (a) \leq  f^{i}_{\beta_k} (a) < f^{1 - i}_{\gamma_k} (a)$}

 whenever $j < k < \cf (\delta)$ and $a \in A \setminus (W_j \cup W_k)$.
\hfill$\square$  

\begin{fact} {\rm(\cite{Weaksat})}  Suppose that there exists an increasing, cofinal sequence $\vec{f} = \langle f_\alpha : \alpha < \kappa \rangle$ in $(\prod A, <_I)$. Then there exists an increasing, cofinal sequence $\vec{g} = \langle g_\alpha : \alpha < \kappa \rangle$ in $(\prod A, <_I)$ such that for any regular cardinal $\sigma$ with $\vert A \vert < \sigma < \sup A$, and any ordinal $\eta$ with $\sigma \leq \eta < \sigma^+$, the following holds.  For any $\beta \in E^\kappa_{\sigma^{+ 3}}$, and any closed unbounded subset $C$ of $\beta$, there is a closed subset $H$ of $C$ of order-type $\eta$ with $H \subseteq G (\vec{g})$.
\end{fact}

\begin{Pro}  Let $\vec{f} = \langle f_\alpha : \alpha < \kappa \rangle$ be an increasing, cofinal sequence in $(\prod A, <_I)$, and $\sigma$ be a regular cardinal with $\vert A \vert < \sigma < \sup A$. Then, setting $J = NS_\kappa \vert (G (\vec{f}) \cap E^\kappa_\sigma)$,  there is either a descending $(J, I_\kappa)$-tower of length ${\frak b}_\kappa$, or an ascending $(J, I_\kappa)$-tower of length $\kappa^+$. 
\end{Pro}

{\bf Proof.}  Set $\rho = \sigma^{+ 3}$, and let $\vec{g} = \langle g_\alpha : \alpha < \kappa \rangle$ be as in the statement of Fact 5.18. 

\medskip

{\bf Claim.}  Let $\gamma \in E^\kappa_{\geq \rho}$. Then $G (\vec{f}) \cap E^\kappa_\sigma \cap \gamma$ is stationary in $\gamma$.

\medskip

{\bf Proof of the claim.} Let $C$ be a closed unbounded subset of $\gamma$. Then we may find $\beta \in E^\kappa_\rho \cap (\gamma + 1)$ such that $C \cap \beta$ is cofinal in $\beta$. There must be a closed subset $H$ of $C \cap \beta$ of order-type $\sigma + 1$ with $H \subseteq G (\vec{g})$. Then $\max H \in C \cap G (\vec{f}) \cap E^\kappa_\sigma$, which completes the proof of the claim.

\medskip

Note that by the claim, $G (\vec{g}) \cap E^\kappa_\sigma$ is stationary in $\kappa$. By Fact 5.17, so is $G (\vec{f}) \cap E^\kappa_\sigma$, and moreover $NS_\kappa \vert (G (\vec{g}) \cap E^\kappa_\sigma) = J$. The desired conclusion is now immediate from Theorem 5.12. 
\hfill$\square$  

\bigskip

\subsection{Robustness of Diamond}

\medskip

Unlike Club (see e.g. \cite{DZ2}), Diamond is remarkably robust, in the sense that a small modification in its definition will often yield an equivalent principle. For a striking example of this, consider the following result which was first established by Primavesi \cite{Prima} for any $J$ of the form $NS_{\omega_1} \vert S$.

\medskip

\begin{Obs} Given a $\kappa$-complete ideal $J$ on $\kappa$ extending $NS_\kappa$, the following are equivalent :
\begin{enumerate}[\rm (i)]
\item $\diamondsuit_\kappa [J]$ holds.  
\item There is $s_\alpha \subseteq \alpha$ for $\alpha < \kappa$ such that $\{\alpha : s_\alpha = C \cap \alpha \} \in J^+$ for all $C \in {\cal C}_\kappa$. 
\end{enumerate}
\end{Obs}

{\bf Proof.} (i) $\rightarrow$ (ii) : Trivial.

(ii) $\rightarrow$ (i) : The proof is a modification of that of Theorem 3.0.10 in \cite{Prima}. Let $\langle s_\alpha : \alpha < \kappa\rangle$ be as in (ii). 

\medskip

{\bf Claim 1.}  $\clubsuit_\kappa [J]$ holds.

\medskip

{\bf Proof of Claim 1.} For $\alpha \in acc (\kappa)$, put $t_\alpha = s_\alpha \setminus acc (s_\alpha)$. Now fix $A \in [\kappa]^\kappa$. Put $C = A \cup acc (A)$, 

\centerline{$D = \{\alpha \in acc (\kappa) : \sup ((C \setminus acc (C)) \cap \alpha) = \alpha\}$}

 and $S = \{\alpha < \kappa : s_\alpha = C \cap \alpha\}$. Clearly for any $\alpha \in D \cap S$, $t_\alpha \subseteq C \setminus acc (C) \subseteq A$, and moreover $\sup t_\alpha \geq \sup ((C \setminus acc (C)) \cap \alpha) = \alpha$, which completes the proof of the claim.

\medskip

{\bf Claim 2.}  $2^{< \kappa} = \kappa$.

\medskip

{\bf Proof of Claim 2.} Let $\tau$ be an infinite cardinal less than $\kappa$. Put $c = acc (\kappa) \cap \tau$, and for any $a \subseteq \tau \setminus acc (\kappa)$, $C_a = c \cup a \cup  (\kappa \setminus \tau)$ and $T_a = \{\alpha < \kappa : s_\alpha = C_a \cap \alpha\}$. Then clearly, $s_\alpha \cap (\tau \setminus acc (\kappa)) = a$ for all $\alpha \in T_a$. It follows that $2^\tau \leq \kappa$, which completes the proof of the claim.

Finally, by Claims 1 and 2 and Observation 2.22, $\diamondsuit_\kappa [J]$ holds.
\hfill$\square$

\medskip

The starred version is established in the same way : 

\medskip

\begin{Obs} Given a $\kappa$-complete ideal $J$ on $\kappa$ extending $NS_\kappa$, the following are equivalent :
\begin{enumerate}[\rm (i)]
\item $\diamondsuit_\kappa^\ast [J]$ holds.  
\item There is $s^{i}_\alpha \subseteq \alpha$ for $i < \alpha < \kappa$ such that $\{\alpha : \exists i < \alpha (s^{i}_\alpha = C \cap \alpha) \} \in J^\ast$ for all $C \in {\cal C}_\kappa$. 
\end{enumerate}
\end{Obs}

\medskip

Not so surprisingly, the situation is different with Club. Suppose for instance that $\kappa = \nu^+$, where $\nu$ is singular, and $J$ is a normal $\kappa^+$-saturated ideal on $\kappa$ (by work of Foreman \cite{Foreman}, this is consistent relative to a huge cardinal). Then by Observation 2.27, $\clubsuit_\kappa^{\rm ev} [J]$ fails. However by Proposition 5.3, we may find $s_\alpha \subseteq \alpha$ with $\sup s_\alpha = \alpha$ for $\alpha \in acc (\kappa)$ such that $\{\alpha : s_\alpha \subseteq C \} \in J^+$ for all $C \in {\cal C}_\kappa$. 

\bigskip

\section{Embarrassing questions}

In this section we attempt to probe the depth of the author's ignorance. Unfortunately, as will shortly be seen, no lower bounds were found.

\begin{Ques} Let $J$ be a normal ideal on $\kappa$ that is not $\kappa^+$-saturated. Does there then exist an ascending $(J, J)$-tower of length $\kappa^+$ ?
\end{Ques}

\begin{Ques} Does $\diamondsuit_\kappa [J]$ imply the existence of a descending (respectively, ascending) $(J, J)$-tower of length $2^\kappa$ ?
\end{Ques}


\medskip

Is it always true that any nontrivial club-like principle for $J$ implies some degree of nonsaturation for $J$ ? Let us consider the following test case. For a $\kappa$-complete ideal $J$ on $\kappa$ extending $NS_\kappa$, and a cardinal $\sigma$ with $2 \leq \sigma \leq \kappa$, let $\clubsuit_\kappa^{\rm ev, \sigma} [J]$ assert the existence of $s_\alpha \subseteq \alpha$ with $\sup s_\alpha = \alpha$ for $\alpha \in acc (\kappa)$ with the property that for any $f : \kappa \rightarrow \sigma$, there is $i < \sigma$ such that

\centerline{$\{\alpha \in acc (\kappa) : \exists \beta < \alpha ((s_\alpha \setminus \beta) \cap f^{- 1} (\{i\}) = \emptyset)\} \in J^+$.}

$\clubsuit_\kappa^{\rm ev, 2} [NS_{\omega_1}]$ is the principle $\clubsuit_{w^2}$ studied in \cite{FSS}. Notice that $\clubsuit_\kappa^{\rm ev, \sigma} [J]$ gets weaker as $\sigma$ increases. 

\medskip

\begin{Obs} Suppose that $J$ is normal and $\clubsuit_\kappa^{\rm ev, \kappa} [J]$ holds. Then $J$ is not prime.
\end{Obs}

{\bf Proof.} Suppose otherwise. Let $s_\alpha \subseteq \alpha$ for $\alpha \in acc (\kappa)$ witness that $\clubsuit_\kappa^{\rm ev, \kappa} [J]$ holds. By a standard argument we may find $S \subseteq \kappa$ such that $T = \{\alpha \in acc (\kappa) : s_\alpha = S \cap \alpha\}$ lies in $J^\ast$. Notice that $\vert S \vert = \kappa$. Select $f : \kappa \rightarrow \kappa$ so that for any $i < \kappa$, $\vert f^{- 1} (\{i\}) \vert = \kappa$, and moreover $\vert  f^{- 1} (\{i\}) \setminus S \vert \leq 1$. Now fix $i < \kappa$. Then clearly $\sup ( f^{- 1} (\{i\}) \cap s_\alpha) = \alpha$ for any $\alpha \in T$ such that $\sup ( f^{- 1} (\{i\}) \cap \alpha) = \alpha$. Contradiction.
\hfill$\square$
 
\begin{Ques} Does $\clubsuit_\kappa^{\rm ev, 2} [J]$ imply that $J$ is not $\kappa$-saturated ? 
\end{Ques}
 
 \bigskip

\section{Turrology / Clubology}

We think that behind each result on nonsaturation, there is a tower (respectively, a club). If we do not see it right away, it does not mean that it is not there, just that more research is needed to find it. It is our opinion that such research is socially useful, as  there should be a tower (respectively, a club) for everyone, not just for higher-ups.

\bigskip

{\bf Acknowledgements}. The author would like to thank Assaf Rinot for fruitful discussions. He is deeply indebted to Moti Gitik for his contribution to this paper.

  \bigskip
\noindent Universit\'e de Caen - CNRS \\
Laboratoire de Math\'ematiques \\
BP 5186 \\
14032 Caen Cedex\\
France\\
Email :  pierre.matet@unicaen.fr\\


\end{document}